\documentclass[11pt,a4paper]{article}
\usepackage[T1]{fontenc}
\usepackage[utf8]{inputenc}
\usepackage{lmodern}
\usepackage[margin=30mm]{geometry}
\usepackage{amsmath,amssymb}
\usepackage{microtype}
\usepackage[hidelinks]{hyperref}
\hypersetup{
  pdftitle={Continuous and compact embeddings for a class of modular sequence spaces},
  pdfauthor={Pierre-A. Vuillermot},
  pdfkeywords={Modular sequence spaces, embeddings, 46B45, 46B50}
}
\begin{document}

\title{Continuous and compact embeddings for a class of modular sequence
spaces}
\author{Pierre-A. Vuillermot$^{\ast ,\ast \ast }$ \\
Center for Mathematical Studies, CEMS.UL, Faculdade de Ci\^{e}ncias,\\
Universidade de Lisboa, 1749-016 Lisboa, Portugal$^{\ast }$\\
Universit\'{e} de Lorraine, CNRS, IECL, F-54000 Nancy,\ France$^{\ast \ast }$
\\
pavuillermot@ciencias.ulisboa.pt}
\date{}
\maketitle

\begin{abstract}
In this article we prove the existence of continuous and compact embeddings
for a class of modular sequence spaces generated by Orlicz functions. The
class under consideration contains sequence spaces of Orlicz-Sobolev type,
including certain special cases that have played an important role in a
number of applications regarding problems arising in non-equilibrium
statistical mechanics, as well as Nakano's class of Lebesgue sequence spaces
with variable exponents. In a few exceptional cases we provide necessary and
sufficient conditions for such embeddings to hold, in many others only
sufficient conditions and in that situation we clearly indicate why these
are not necessary by means of additional criteria and specific examples. Our
results are complementary to those already available in the literature. We
illustrate them with several concrete examples.


\end{abstract}

\section{Introduction and outline}

Continuous and compact embeddings involving functional spaces of various
kinds play an important role in Mathematical Analysis and its applications
(see, e.g., \cite{Edmundsevans}).\ This includes the sequence spaces
recently used in the investigation of certain questions in non-equilibrium
statistical mechanics (see, e.g., \cite{Vuillermot} and its Appendix A).
Those spaces are in fact very particular examples of the class of modular
spaces defined and analyzed in the forthcoming sections. It is also worth
noting that some embeddings involving them were already generalized in \cite%
{Vuillermotbis} and \cite{Vuillermotter}. In this article we go one step
further by introducing a class of functional spaces generated by sequences
of Orlicz functions. Accordingly, we organize the remaining part of this
paper in the following way: In Section 2 we define the modular class we are
interested in and the associated modular spaces. We discuss various features
of them with an almost exclusive emphasis on those properties that are used
in Section 3. All the properties listed there are known (see, e.g., standard
references on the subject such as \cite{Krasnoselskiirut}-\cite{Raorenbis}),
but we make a point of either sketching the proof of some of them in order
to make this article self-contained and thereby more accessible to the
interested reader, or else provide the relevant bibliography. We also note
that the modular class we introduce contains as a particular case Nakano's
class of Lebesgue sequence spaces with variable exponents introduced in \cite%
{Nakano}, whose embedding properties and related matters have been the
subject of many investigations over the years up to very recently (see,
e.g., \cite{Langnek} and its references). In Section 3 we prove several
continuous and compact embeddings involving the modular spaces introduced in
Section 2. In some exceptional cases the conditions we impose for such
embeddings to hold are necessary and sufficient, in many others only
sufficient and in that situation we provide examples and other criteria
showing why the conditions are not necessary. We also show that some of the
theorems we state encompass earlier results regarding weighted
Orlicz-Sobolev type spaces. To the best of our knowledge, our results are
new and complementary to those already available in the literature. In this
regard we refer the reader to Sections 8 and 9 in Chapter II of \cite%
{Musielak} and to the numerous references listed there. Whereas our results
focus exclusively on modular spaces of sequences, those sections contain a
plethora of statements valid for modular spaces defined on general measure
spaces, including in particular a generalization of Kolmogorov's and Riesz's
compactness criteria based on properties of Steklov functions (for an
alternative exposition of such criteria in the simpler context of Orlicz
spaces, see Section 11 in Chapter II of \cite{Krasnoselskiirut}). Other
criteria for embeddings can also be found in \cite{Diharu}, \cite{Raoren}
and \cite{Raorenbis}. In particular, as far as earlier results regarding
weak sequential compactness in Orlicz sequence spaces are concerned, we
refer to Section 10.3 in Chapter X of \cite{Raorenbis}. More recently,
regarding modular spaces defined with respect to discrete or counting
measures, we also mention \cite{banomaak} whose authors developed relative
compactness criteria in a fairly general context with a special emphasis on
an analysis of the discrete Hardy operator.

\section{Variations on known facts}

The spaces we alluded to in the introduction are generated by sequences $%
\left( \Phi _{\mathsf{m}}\right) _{\mathsf{m}\in \mathbb{Z}}$ of continuous
functions $\Phi _{\mathsf{m}}:\left[ 0,\infty \right) \mapsto \left[
0,\infty \right) $ that are convex and satisfy $\Phi _{\mathsf{m}}\left(
0\right) =0,\Phi _{\mathsf{m}}\left( t\right) >0$ for every $t>0$. It
follows immediately from these properties that all $\Phi _{\mathsf{m}}$ are
strictly increasing and that $\lim_{t\rightarrow \infty }\Phi _{\mathsf{m}%
}\left( t\right) =\infty $. We refer to them as \textit{Orlicz functions }in
the sequel and simply write $\hat{\Phi}:=$ $\left( \Phi _{\mathsf{m}}\right)
_{\mathsf{m}\in \mathbb{Z}}$ or $\left( \Phi _{\mathsf{m}}\right) $ for the
corresponding sequence when no confusion arises. It is worth noting that
this notion is all we need at the moment. Two crucial additional
restrictions will be imposed later. We begin with the following:

\bigskip

\textbf{Definition 1.} Let $\left( \Phi _{\mathsf{m}}\right) $ be a sequence
of Orlicz functions. We write $\mathsf{c}_{\hat{\Phi}}$ for the \textit{%
modular sequence class} consisting of all complex sequences $\mathsf{p}=$ $%
\left( p_{\mathsf{m}}\right) $ satisfying%
\begin{equation}
\sum_{\mathsf{m}\in \mathbb{Z}}\Phi _{\mathsf{m}}\left( \left\vert p_{%
\mathsf{m}}\right\vert \right) <\infty .  \label{modularclass}
\end{equation}

It follows easily from the monotonicity and the convexity properties of the $%
\Phi _{\mathsf{m}}$ that $\mathsf{c}_{\hat{\Phi}}$ is a convex and balanced
set in the sense that for all $\mathsf{p,q\in c}_{\hat{\Phi}}$ we have $%
\lambda \mathsf{p+\mu q\in c}_{\hat{\Phi}}$ whenever $\lambda ,\mu \in 
\mathbb{C}$ satisfy $\left\vert \lambda \right\vert +\left\vert \mu
\right\vert \leq 1.$ In particular, if $\mu =0$ and $\lambda \in \left[ -1,1%
\right] $ then $\lambda \mathsf{p\in c}_{\hat{\Phi}}$. However, $\mathsf{c}_{%
\hat{\Phi}}$ is not a vector space in general. For instance, let us consider
the Orlicz function 
\begin{equation}
\Phi \left( t\right) =\left\{ 
\begin{array}{c}
0\text{ \ if \ }t=0 \\ 
\\ 
\exp \left[ t-\frac{1}{t}\right] \text{ \ if \ }t>0%
\end{array}%
\right.  \label{orliczfunction}
\end{equation}%
and the constant sequence $\left( \Phi _{\mathsf{m}}\right) $ where $\Phi _{%
\mathsf{m}}=\Phi $ for every $\mathsf{m}\in \mathbb{Z}$. Let $\mathsf{p}=$ $%
\left( p_{\mathsf{m}}\right) $ be given by%
\begin{equation}
p_{\mathsf{m}}=\left\{ 
\begin{array}{c}
0\text{ \ if \ }\mathsf{m}\in \mathbb{Z}\cap \left( -\infty ,0\right] , \\ 
\\ 
\frac{1}{\ln \left( \mathsf{m}^{2}+1\right) }\text{ \ if }\mathsf{m}\in 
\mathbb{N}.%
\end{array}%
\right.  \label{sequencebis}
\end{equation}%
Then we have 
\begin{equation*}
\sum_{\mathsf{m}\in \mathbb{Z}}\Phi \left( \left\vert p_{\mathsf{m}%
}\right\vert \right) =\sum_{\mathsf{m}=1}^{\infty }\frac{\exp \left[ \frac{1%
}{\ln \left( \mathsf{m}^{2}+1\right) }\right] }{\mathsf{m}^{2}+1}<\infty ,
\end{equation*}%
whereas 
\begin{equation*}
\sum_{\mathsf{m}\in \mathbb{Z}}\Phi \left( 2\left\vert p_{\mathsf{m}%
}\right\vert \right) =\sum_{\mathsf{m}=1}^{\infty }\frac{\exp \left[ \frac{2%
}{\ln \left( \mathsf{m}^{2}+1\right) }\right] }{\left( \mathsf{m}%
^{2}+1\right) ^{\frac{1}{2}}}=\infty .
\end{equation*}%
\newline
Therefore $\mathsf{p\in c}_{\hat{\Phi}}$ but $2\mathsf{p\notin c}_{\hat{\Phi}%
}$. Necessary and sufficient conditions that make $\mathsf{c}_{\hat{\Phi}}$
a linear space will be recalled below.

For a given sequence of Orlicz functions $\left( \Phi _{\mathsf{m}}\right) $
we now define the \textit{large modular sequence space} $l_{\hat{\Phi}}$ as
the set of all sequences $\mathsf{p}=$ $\left( p_{\mathsf{m}}\right) $ for
which there exists a finite $\rho ^{\ast }>0$ generally depending on $%
\mathsf{p}$ such that%
\begin{equation}
\sum_{\mathsf{m}\in \mathbb{Z}}\Phi _{\mathsf{m}}\left( \frac{\left\vert p_{%
\mathsf{m}}\right\vert }{\rho ^{\ast }}\right) <\infty .
\label{convergencesexto}
\end{equation}%
We give a proof of the first part of the following result as it provides
another possible definition of $l_{\hat{\Phi}}$:

\bigskip

\textbf{Theorem 1. }\textit{The set }$l_{\hat{\Phi}}$\textit{\ is a complex
Banach space relative to the usual algebraic operations on sequences and
Luxemburg's gauge norm}

\begin{equation}
\left\Vert \mathsf{p}\right\Vert _{\hat{\Phi}}:=\inf \left\{ \rho >0:\sum_{%
\mathsf{m}\in \mathbb{Z}}\Phi _{\mathsf{m}}\left( \frac{\left\vert p_{%
\mathsf{m}}\right\vert }{\rho }\right) \leq 1\right\} .  \label{norm}
\end{equation}%
\textit{Moreover, we have} 
\begin{equation*}
l_{\hat{\Phi}}=\operatorname{span}\mathsf{c}_{\hat{\Phi}}.
\end{equation*}

\bigskip

\textbf{Proof.} For any $\mathsf{p}$ $\in l_{\hat{\Phi}}$ we define the
function $F_{\mathsf{p}}:\left[ \rho ^{\ast },\infty \right) \mapsto \left[
0,\infty \right) $ by 
\begin{equation}
F_{\mathsf{p}}\left( \rho \right) :=\sum_{\mathsf{m}\in \mathbb{Z}}\Phi _{%
\mathsf{m}}\left( \frac{\left\vert p_{\mathsf{m}}\right\vert }{\rho }\right)
.  \label{function}
\end{equation}%
If $\mathsf{p=0}$ we have $F_{\mathsf{p}}=0$ identically. If $\mathsf{p\neq 0%
}$ then $F_{\mathsf{p}}$ remains finite and is strictly decreasing as a
consequence of (\ref{convergencesexto}), for each $\Phi _{\mathsf{m}}$ is
strictly increasing. We wish to prove that%
\begin{equation}
\lim_{\rho \rightarrow \infty }F_{\mathsf{p}}\left( \rho \right) =0.
\label{convergencetozero}
\end{equation}%
Let $\left( \rho _{\mathsf{n}}\right) $ be any sequence such that $\rho _{%
\mathsf{n}}\in \left[ \rho ^{\ast },\infty \right) $ for every $\mathsf{n}%
\in \mathbb{N}$ and satisfying $\rho _{\mathsf{n}}\rightarrow \infty $ as $%
\mathsf{n\rightarrow \infty }$. Then for every $\mathsf{m}\in \mathbb{Z}$ we
get%
\begin{equation*}
\Phi _{\mathsf{m}}\left( \frac{\left\vert p_{\mathsf{m}}\right\vert }{\rho _{%
\mathsf{n}}}\right) \rightarrow 0
\end{equation*}%
as $\mathsf{n\rightarrow \infty }$, and%
\begin{equation*}
\Phi _{\mathsf{m}}\left( \frac{\left\vert p_{\mathsf{m}}\right\vert }{\rho _{%
\mathsf{n}}}\right) \leq \Phi _{\mathsf{m}}\left( \frac{\left\vert p_{%
\mathsf{m}}\right\vert }{\rho ^{\ast }}\right)
\end{equation*}%
uniformly in $\mathsf{n}$ on the right-hand side. Keeping in mind that (\ref%
{convergencesexto}) holds\bigskip , we then obtain%
\begin{equation*}
\lim_{\mathsf{n}\rightarrow \infty }F_{\mathsf{p}}\left( \rho _{\mathsf{n}%
}\right) =0
\end{equation*}%
from dominated convergence. Therefore (\ref{convergencetozero}) holds and we
have $F_{\mathsf{p}}\left( \rho \right) \leq 1$ for all sufficiently large $%
\rho $, so that the set on the right-hand side of (\ref{norm}) is non-empty
and thus $\left\Vert \mathsf{p}\right\Vert _{\hat{\Phi}}<\infty $.

Conversely, if $\left\Vert \mathsf{p}\right\Vert _{\hat{\Phi}}<\infty $ then
the set in question is non-empty so that there exists \textit{a fortiori} at
least one $\rho ^{\ast }>0$ such that (\ref{convergencesexto}) holds.
Consequently, we may identify $l_{\hat{\Phi}}$ with the set of \textit{all}
sequences such that $\left\Vert \mathsf{p}\right\Vert _{\hat{\Phi}}<\infty $%
. But it is known that (\ref{norm}) has all the attributes of a norm, which
implies that $l_{\hat{\Phi}}$ is a normed vector space with respect to the
usual algebraic operations on sequences. The fact that $l_{\hat{\Phi}}$ is
complete with respect to (\ref{norm}) follows from standard arguments (see,
e.g., \cite{Krasnoselskiirut}-\cite{Raorenbis}).

Let us now consider $\mathsf{p}$ $\in \mathsf{c}_{\hat{\Phi}}$. A reasoning
similar to the one above then shows that $\left\Vert \mathsf{p}\right\Vert _{%
\hat{\Phi}}<\infty $ as well, so that $\mathsf{c}_{\hat{\Phi}}\subseteq l_{%
\hat{\Phi}}$ and hence%
\begin{equation*}
\operatorname{span}\mathsf{c}_{\hat{\Phi}}\subseteq l_{\hat{\Phi}}\text{.}
\end{equation*}%
Finally, if $\mathsf{p}$ $\in l_{\hat{\Phi}}$ there exists $\rho ^{\ast }>0$
such that (\ref{convergencesexto}) holds, and so with $\mathsf{q:=}\frac{%
\mathsf{p}}{\rho ^{\ast }}$ we have $\mathsf{q\in c}_{\hat{\Phi}}$ and
thereby $\mathsf{p=}\rho ^{\ast }\mathsf{q\in }\operatorname{span}\mathsf{c}_{%
\hat{\Phi}}$. Consequently we have%
\begin{equation*}
\ \ l_{\hat{\Phi}}\subseteq \operatorname{span}\mathsf{c}_{\hat{\Phi}}.
\end{equation*}%
$\blacksquare $

\bigskip

Aside from $l_{\hat{\Phi}}$ it is also possible to define the \textit{small
modular sequence space} $\mathsf{s}_{\hat{\Phi}}$ consisting of all complex
sequences $\mathsf{p}=$ $\left( p_{\mathsf{m}}\right) $ satisfying%
\begin{equation*}
\sum_{\mathsf{m}\in \mathbb{Z}}\Phi _{\mathsf{m}}\left( \frac{\left\vert p_{%
\mathsf{m}}\right\vert }{\rho }\right) <\infty
\end{equation*}%
for \textit{every} $\rho >0$. From this and from the monotonicity and
convexity properties of the $\Phi _{\mathsf{m}}$, it follows that $\mathsf{s}%
_{\hat{\Phi}}$ is a complex vector space with respect to the usual algebraic
operations. Furthermore we have%
\begin{equation*}
\mathsf{s}_{\hat{\Phi}}\subseteq \mathsf{c}_{\hat{\Phi}}\subseteq l_{\hat{%
\Phi}},
\end{equation*}%
these inclusions being strict in general as shown by the above example of
the constant sequence $\left( \Phi _{\mathsf{m}}\right) $ generated by (\ref%
{orliczfunction}). More to the point, from (\ref{orliczfunction}), (\ref%
{sequencebis}) and $\rho >0$ we have%
\begin{equation*}
\sum_{\mathsf{m}\in \mathbb{Z}}\Phi _{\mathsf{m}}\left( \frac{\left\vert p_{%
\mathsf{m}}\right\vert }{\rho }\right) =\sum_{\mathsf{m}=1}^{\infty }\frac{%
\exp \left[ \frac{1}{\rho \ln \left( \mathsf{m}^{2}+1\right) }\right] }{%
\left( \mathsf{m}^{2}+1\right) ^{\rho }}<\infty
\end{equation*}%
if, and only if, $\rho \in \left( \frac{1}{2},\infty \right) $, so that $%
\mathsf{p}\in \mathsf{c}_{\hat{\Phi}}$ but $\mathsf{p}\notin \mathsf{s}_{%
\hat{\Phi}}$.

In the following result we state additional properties regarding $\mathsf{s}%
_{\hat{\Phi}}$ which will be quite relevant in the following section, and we
propose a short sketch of the proof:

\bigskip

\textbf{Theorem 2. }\textit{Let }$\left( \Phi _{\mathsf{m}}\right) $\textit{%
\ be a sequence of Orlicz functions. Then the linear space }$\mathsf{s}_{%
\hat{\Phi}}$ \textit{is a complex Banach space with respect to (\ref{norm}),
and is the largest subspace contained in the modular sequence class }$%
\mathsf{c}_{\hat{\Phi}}$\textit{. Moreover, there exists a Schauder basis in 
}$\mathsf{s}_{\hat{\Phi}}$\textit{.}

\bigskip

\textbf{Proof. }The fact that $\mathsf{s}_{\hat{\Phi}}$ is a closed subspace
of $l_{\hat{\Phi}}$ follows from standard arguments (see, e.g. \cite%
{Krasnoselskiirut}-\cite{Raorenbis} and \cite{Vuillermotter}). Furthermore,
let $\mathsf{S}$ be a linear subspace of $\mathsf{c}_{\hat{\Phi}}$, let $%
\rho >0$ be arbitrary and let us choose $\lambda \in \mathbb{C}$ such that $%
\left\vert \lambda \right\vert ^{-1}=\rho $. For any $\mathsf{p\in S}$ we
then have $\lambda \mathsf{p}\in \mathsf{S}\subseteq $ $\mathsf{c}_{\hat{\Phi%
}}$ and thereby%
\begin{equation*}
\sum_{\mathsf{m}\in \mathbb{Z}}\Phi _{\mathsf{m}}\left( \frac{\left\vert p_{%
\mathsf{m}}\right\vert }{\rho }\right) =\sum_{\mathsf{m}\in \mathbb{Z}}\Phi
_{\mathsf{m}}\left( \left\vert \left( \lambda \mathsf{p}\right) _{\mathsf{m}%
}\right\vert \right) <\infty ,
\end{equation*}%
which shows that $\mathsf{S}\subseteq $ $\mathsf{s}_{\hat{\Phi}}$. Finally,
it is the set of sequences $\mathsf{e}_{\mathsf{m}}$ defined by $\left( 
\mathsf{e}_{\mathsf{m}}\right) _{\mathsf{n}}=\delta _{\mathsf{m,n}}$ for all 
$\mathsf{m,n}\in \mathbb{Z}$ that constitutes a Schauder basis of $\mathsf{s}%
_{\hat{\Phi}}$. Indeed, it is easily verified that $\mathsf{e}_{\mathsf{m}%
}\in \mathsf{s}_{\hat{\Phi}}$ for each $\mathsf{m}$ so that for every $%
\mathsf{M}\in \mathbb{N}$ we may define%
\begin{equation}
\mathsf{p}_{\mathsf{M}}:\mathsf{=}\sum_{\left\vert \mathsf{m}\right\vert
\leq \mathsf{M}}p_{\mathsf{m}}\mathsf{e}_{\mathsf{m}}  \label{cauchysequence}
\end{equation}%
as an element of that space. One then shows that (\ref{cauchysequence})
generates a Cauchy sequence therein. This means in particular that every $%
\mathsf{p}\in \mathsf{s}_{\hat{\Phi}}$ may be expanded in a unique fashion
as 
\begin{equation*}
\mathsf{p=}\sum_{\mathsf{m}\in \mathbb{Z}}p_{\mathsf{m}}\mathsf{e}_{\mathsf{m%
}}\text{,}
\end{equation*}%
a norm-convergent series in $\mathsf{s}_{\hat{\Phi}}$. \ \ $\blacksquare $

\bigskip

\textsc{Remarks.} (1) Whereas $\left( \mathsf{e}_{\mathsf{m}}\right) $
provides a basis for $\mathsf{s}_{\hat{\Phi}}$, it is worth keeping in mind
that there is no Schauder basis in $l_{\hat{\Phi}}$, nor is this space
separable in general unless of course $\mathsf{s}_{\hat{\Phi}}=l_{\hat{\Phi}%
} $. This raises the important question of characterizing those sequences of
Orlicz functions $\left( \Phi _{\mathsf{m}}\right) $ that lead to such an
identification, for the existence of a Schauder basis will be central to the
developement of relative compactness arguments in the next section. We shall
come back to this problem shortly.

(2) According to the definition of the gauge norm (\ref{norm}), there exists
a minimizing sequence of positive numbers $\left( \rho _{\mathsf{n}}\right) $
satisfying%
\begin{equation*}
\sum_{\mathsf{m}\in \mathbb{Z}}\Phi _{\mathsf{m}}\left( \frac{\left\vert p_{%
\mathsf{m}}\right\vert }{\rho _{\mathsf{n}}}\right) \leq 1
\end{equation*}%
for each $\mathsf{n}\in \mathbb{N}$ and $\rho _{\mathsf{n}}\rightarrow
\left\Vert \mathsf{p}\right\Vert _{\hat{\Phi}}$ as $\mathsf{n\rightarrow
\infty }$, which we may assume to be monotone decreasing. Therefore, we
infer from monotone convergence that%
\begin{equation}
\sum_{\mathsf{m}\in \mathbb{Z}}\Phi _{\mathsf{m}}\left( \frac{\left\vert p_{%
\mathsf{m}}\right\vert }{\left\Vert \mathsf{p}\right\Vert _{\hat{\Phi}}}%
\right) \leq 1  \label{importantrelation}
\end{equation}%
for every $\mathsf{p}\in l_{\hat{\Phi}}$ such that $\left\Vert \mathsf{p}%
\right\Vert _{\hat{\Phi}}>0$. Moreover, if $\left\Vert \mathsf{p}\right\Vert
_{\hat{\Phi}}\leq 1$ we have 
\begin{equation}
\sum_{\mathsf{m}\in \mathbb{Z}}\Phi _{\mathsf{m}}\left( \left\vert p_{%
\mathsf{m}}\right\vert \right) \leq \left\Vert \mathsf{p}\right\Vert _{\hat{%
\Phi}}  \label{importantrelationbis}
\end{equation}%
which is obvious if $\mathsf{p=0}$, while if $\mathsf{p\neq 0}$ it is a
consequence of the convexity of $\Phi _{\mathsf{m}}$, the fact that $\Phi _{%
\mathsf{m}}$ vanishes at the origin together with Relation (\ref%
{importantrelation}). Both (\ref{importantrelation}) and (\ref%
{importantrelationbis}) play a crucial role in Section 3.

\bigskip

Many properties of the modular spaces defined above including the existence
of certain embeddings are intimately tied up with the following conditions
that govern the behavior of the Orlicz functions either in a small
neighborhood of the origin or globally (see, e.g., Section d. of Chapter 4
in \cite{Lindenstrausstza} for related formulations):

\bigskip

\textbf{Definition 2. }(1) We say the sequence of Orlicz functions $\left(
\Phi _{\mathsf{m}}\right) $ satisfies the \textit{uniform }$\Delta _{2}$%
\textit{-condition at the origin} if there exist finite constants $\eta
>0,t_{0}>0$, both independent of $\mathsf{m,}$ and $\mathsf{M}\in \mathbb{N}$
such that the inequality%
\begin{equation}
\Phi _{\mathsf{m}}\left( 2t\right) \leq \eta \Phi _{\mathsf{m}}\left(
t\right)  \label{inequalityquinto}
\end{equation}%
holds for each $t\in \left[ 0,t_{0}\right] $ and every $\mathsf{m}\in 
\mathbb{Z}$ with $\left\vert \mathsf{m}\right\vert \geq \mathsf{M}.$

(2) We say the sequence of Orlicz functions $\left( \Phi _{\mathsf{m}%
}\right) $ satisfies the \textit{uniform global }$\Delta _{2}$\textit{%
-condition} if there exists a finite constant $\eta >0$ independent of $%
\mathsf{m}$ and $\mathsf{M}\in \mathbb{N}$ such that the inequality%
\begin{equation}
\Phi _{\mathsf{m}}\left( 2t\right) \leq \eta \Phi _{\mathsf{m}}\left(
t\right)  \label{globalcondition}
\end{equation}%
holds for each real $t\geq 0$ and every $\mathsf{m}\in \mathbb{Z}$ with $%
\left\vert \mathsf{m}\right\vert \geq \mathsf{M}.$

\bigskip

\textsc{Remark.} It follows from (1) of the above definition that there
exist an $\mathsf{M}\in \mathbb{N}$, and for every $\theta \in \left(
0,\infty \right) $ finite constants $\eta _{\theta }>0,t_{\theta }>0$
independent of $\mathsf{M}$, such that%
\begin{equation}
\Phi _{\mathsf{m}}\left( \theta t\right) \leq \eta _{\theta }\Phi _{\mathsf{m%
}}\left( t\right)  \label{inequalitysexto}
\end{equation}%
for each $t\in \left[ 0,t_{\theta }\right] $ and every $\mathsf{m}\in 
\mathbb{Z}$ with $\left\vert \mathsf{m}\right\vert \geq \mathsf{M}.$ Indeed,
inequality (\ref{inequalityquinto}) first extends by iteration to 
\begin{equation*}
\Phi _{\mathsf{m}}\left( 2^{\mathsf{n}}t\right) \leq \eta ^{\mathsf{n}}\Phi
_{\mathsf{m}}\left( t\right)
\end{equation*}%
for each $t\in \left[ 0,2^{1-\mathsf{n}}t_{0}\right] $ and every $\mathsf{n}%
\in \mathbb{N}$, and then to (\ref{inequalitysexto}) since with $\mathsf{n}%
_{\theta }\geq \log _{2}\theta $ we have $\theta \leq 2^{\mathsf{n}_{\theta
}}$. In other words, (\ref{inequalitysexto}) holds for every $\theta \in
\left( 0,\infty \right) $ with constants depending on $\theta $ in general.
The converse statement, namely, that (\ref{inequalitysexto}) implies (\ref%
{inequalityquinto}), is obvious but it is in the form (\ref{inequalitysexto}%
) that the $\Delta _{2}$\textit{-}condition will be used. Similarly, (\ref%
{globalcondition}) holds if, and only if, the inequality%
\begin{equation}
\Phi _{\mathsf{m}}\left( \theta t\right) \leq \eta _{\theta }\Phi _{\mathsf{m%
}}\left( t\right)  \label{inequalitysextobis}
\end{equation}%
holds for each $\theta \in \left( 0,\infty \right) $ and every $t\geq 0$ for
all $\mathsf{m}\in \mathbb{Z}$ with $\left\vert \mathsf{m}\right\vert \geq 
\mathsf{M}.$

\bigskip

A typical application is the following result:

\bigskip

\textbf{Theorem 3.} \textit{We have }%
\begin{equation}
\mathsf{s}_{\hat{\Phi}}=\mathsf{c}_{\hat{\Phi}}=l_{\hat{\Phi}}
\label{banachspace}
\end{equation}%
\textit{if, and only if, the sequence }$\left( \Phi _{\mathsf{m}}\right) $ 
\textit{satisfies the uniform }$\Delta _{2}$\textit{-condition at the
origin. In this case the Banach space in (\ref{banachspace}) consists of all
complex sequences }$\mathsf{(}p\mathsf{_{\mathsf{m}})}$\textit{\ satisfying}%
\begin{equation*}
\sum_{\mathsf{m}\in \mathbb{Z}}\Phi _{\mathsf{m}}\left( \left\vert p_{%
\mathsf{m}}\right\vert \right) <\infty
\end{equation*}%
\textit{and has }$\left( \mathsf{e}_{\mathsf{m}}\right) $ \textit{as a
Schauder basis.}

\bigskip

\textbf{Proof. }The first statement follows from the standard arguments
developed in \cite{Krasnoselskiirut}-\cite{Raorenbis} and \cite%
{Vuillermotter}. The proof of the second statement is an immediate
consequence of (\ref{modularclass}) and Theorem 2. \ \ $\blacksquare $

\bigskip

\textsc{Remarks.} (1) It is worth recalling that the uniform $\Delta _{2}$%
-condition at the origin does not prevent Orlicz functions from growing
exponentially rapidly when the independent variable becomes large. Thus, the
functions%
\begin{equation*}
\Phi _{\mathsf{m}}\left( t\right) =\frac{\exp \left[ t\right] -t-1}{%
\left\vert \mathsf{m}\right\vert +1}
\end{equation*}%
do satisfy the condition, whereas (\ref{orliczfunction}) does not. On the
other hand, the functions%
\begin{equation*}
\Phi _{\mathsf{m}}\left( t\right) =t^{\alpha _{\mathsf{m}}}
\end{equation*}%
with $\alpha _{\mathsf{m}}\in \left[ 1,\infty \right) $ satisfy the uniform
global $\Delta _{2}$-condition when the sequence $\left( \alpha _{\mathsf{m}%
}\right) $ is bounded. The corresponding modular spaces are then Nakano's we
briefly alluded to in the introduction. While their continuous embeddings
and related properties are thoroughly discussed in \cite{Langnek}, we shall
give a simple criterion for the relative compactness of their bounded sets
in the next section.

(2) The functions%
\begin{equation*}
\left\{ 
\begin{array}{c}
\Phi \left( t\right) =\exp \left[ t\right] -1, \\ 
\\ 
\Phi \left( t\right) =t%
\end{array}%
\right.
\end{equation*}%
are both Orlicz functions in the sense of the definition given at the
beginning of this section, but their respective Legendre transform computed
from%
\begin{equation}
\Phi ^{\ast }\left( s\right) =\sup_{t\in \left[ 0,\infty \right) }\left(
st-\Phi \left( t\right) \right) \text{, \ \ }s\geq 0\text{,}
\label{legendretransform}
\end{equation}%
are not. Indeed, in the first case we get 
\begin{equation*}
\Phi ^{\ast }\left( s\right) =s\ln s-s+1,
\end{equation*}%
while in the second case we have%
\begin{equation*}
\Phi ^{\ast }\left( s\right) =\left\{ 
\begin{array}{c}
0\text{ \ \ if }0\leq s\leq 1, \\ 
\\ 
\infty \text{ \ if }s\in \left( 1,\infty \right) .%
\end{array}%
\right.
\end{equation*}%
However, we shall see in one instance of Section 3 that it is sometimes
quite advantageous to have some kind of criterion that characterizes the
weak convergence to zero of the Schauder basis $\left( \mathsf{e}_{\mathsf{m}%
}\right) $ in $\mathsf{s}_{\hat{\Phi}}$, and this calls for a smaller class
of functions that does remain invariant under (\ref{legendretransform}).
What we have in mind is the following statement, whose proof we sketch for
the sake of completeness:

\bigskip

\textbf{Proposition 1. }\textit{Let }$\left( \Phi _{\mathsf{m}}\right) $%
\textit{\ be a sequence of Orlicz functions, and let us assume that each
right-derivative }$\Phi _{r,\mathsf{m}}$ \textit{satisfies }$\Phi _{r,%
\mathsf{m}}\left( 0\right) =0$ \textit{and }$\lim_{t\rightarrow \infty }\Phi
_{r,\mathsf{m}}\left( t\right) =\infty $. \textit{Then the sequence of
Legendre transforms} $\left( \Phi _{\mathsf{m}}^{\ast }\right) $ \textit{%
also consists of Orlicz functions satisfying these additional properties.
Under these conditions, the sequence of basis vectors }$\left( \mathsf{e}_{%
\mathsf{m}}\right) $ \textit{converges weakly to zero in} $\mathsf{s}_{\hat{%
\Phi}}$ \textit{if, and only if, }%
\begin{equation}
\liminf_{\left\vert \mathsf{m}\right\vert \rightarrow \infty }\Phi _{\mathsf{%
m}}^{\ast }\left( t\right) >0  \label{constraintbis}
\end{equation}%
\textit{for every} $t>0$.

\bigskip

\textbf{Proof. }The first part of the statement is well known, and is at the
basis of the duality theory for Orlicz and modular sequence spaces (see,
e.g., Sections b. and d. of Chapter 4 in \cite{Lindenstrausstza}, and the
references therein).

As for the statement regarding weak convergence, we begin by observing that
the sequence of basis vectors $\left( \mathsf{e}_{\mathsf{m}}\right) $
converges weakly to zero in $s_{\hat{\Phi}}$ if, and only if,%
\begin{equation}
l_{\hat{\Phi}^{\ast }}\subset c_{0},  \label{inclusionbis}
\end{equation}%
where $\hat{\Phi}^{\ast }:=$ $\left( \Phi _{\mathsf{m}}^{\ast }\right) _{%
\mathsf{m}\in \mathbb{Z}}$ and where $c_{0}$ stands for the space of all
complex sequences converging to zero. Indeed, as the strong dual of $s_{\hat{%
\Phi}}$ indentifies with $l_{\hat{\Phi}^{\ast }}$, every linear bounded
functional on $s_{\hat{\Phi}}$ is generated by a $\mathsf{q}\in l_{\hat{\Phi}%
^{\ast }}$ and is of the form%
\begin{equation*}
\varphi _{\mathsf{q}}\left( \mathsf{p}\right) =\displaystyle\sum\limits_{\mathsf{m}\in 
\mathbb{Z}}\bar{q}_{\mathsf{m}}p_{\mathsf{m}}
\end{equation*}%
for each $\mathsf{p}\in s_{\hat{\Phi}}$, the preceding series being
absolutely convergent. Therefore, $\varphi _{\mathsf{q}}\left( \mathsf{e}_{%
\mathsf{m}}\right) =\bar{q}_{\mathsf{m}}$ for all $\mathsf{q}\in l_{\hat{\Phi%
}^{\ast }}$ so that the desired equivalence follows.

It remains to show that (\ref{constraintbis}) is equivalent to (\ref%
{inclusionbis}). Thus, let us assume that (\ref{constraintbis}) does not
hold, that is, there exists a finite $t^{\ast }>0$ such that%
\begin{equation*}
\liminf_{\left\vert \mathsf{m}\right\vert \rightarrow \infty }\Phi _{\mathsf{%
m}}^{\ast }\left( t^{\ast }\right) =0.
\end{equation*}%
Then there is a subsequence indexed by $\left( \mathsf{m}_{\mathsf{n}%
}\right) _{\mathsf{n}\in \mathbb{N}}$ such that%
\begin{equation}
\Phi _{\mathsf{m}_{\mathsf{n}}}^{\ast }\left( t^{\ast }\right) \leq 2^{-%
\mathsf{n}}  \label{majoration}
\end{equation}%
for all sufficiently large $\mathsf{n}$. Therefore, defining $\mathsf{p}%
=\left( p_{\mathsf{m}}\right) $ by%
\begin{equation*}
p_{\mathsf{m}}=\left\{ 
\begin{array}{c}
t^{\ast }\text{ \ if }\mathsf{m}=\mathsf{m}_{\mathsf{n}} \\ 
\\ 
0\text{ \ if }\mathsf{m}\neq \mathsf{m}_{\mathsf{n}}%
\end{array}%
\right.
\end{equation*}%
we see that $\mathsf{p}\notin c_{0}$. But on the other hand we have%
\begin{equation*}
\sum_{\mathsf{m}\in \mathbb{Z}}\Phi _{\mathsf{m}}^{\ast }\left( \left\vert
p_{\mathsf{m}}\right\vert \right) =\sum_{\mathsf{n}=1}^{\infty }\Phi _{%
\mathsf{m}_{\mathsf{n}}}^{\ast }\left( t^{\ast }\right) <\infty
\end{equation*}%
by virtue of (\ref{majoration}) so that $\mathsf{p}\in l_{\hat{\Phi}^{\ast
}} $, thereby implying that (\ref{inclusionbis}) does not hold. Consequently
(\ref{inclusionbis}) implies (\ref{constraintbis}). The converse inclusion
follows from similar arguments. \ \ $\blacksquare $

(3) Let $\Phi $ be an Orlicz function, let $w:=\left( w_{\mathsf{m}}\right) $
be a sequence of positive numbers and let $k\in \mathbb{R}$. We consider the
very particular sequence of functions $\left( \Phi _{\mathsf{m}}\right) $
given by%
\begin{equation}
\Phi _{\mathsf{m}}\left( t\right) :=w_{\mathsf{m}}\left( 1+\Phi \left(
\left\vert \mathsf{m}\right\vert \right) \right) ^{k}\Phi \left( t\right) ,
\label{exempleparticulier}
\end{equation}%
which gives rise to the modular class $\mathsf{c}_{k,\Phi ,w}$ consisting of
all complex sequences $\left( p_{\mathsf{m}}\right) $ satisfying%
\begin{equation*}
\sum_{\mathsf{m}\in \mathbb{Z}}w_{\mathsf{m}}\left( 1+\Phi \left( \left\vert 
\mathsf{m}\right\vert \right) \right) ^{k}\Phi \left( \left\vert p_{\mathsf{m%
}}\right\vert \right) <\infty .
\end{equation*}%
The corresponding modular spaces $\mathsf{s}_{\hat{\Phi}}$, $l_{\hat{\Phi}}$
are then written $\mathsf{s}_{k,\Phi ,w}$, $l_{k,\Phi ,w}$ for $k\neq 0$,
and simply $\mathsf{s}_{\Phi ,w}$, $l_{\Phi ,w}$ when $k=0$, respectively,
and already appear in \cite{Vuillermotter} as \textit{Orlicz-Sobolev type
spaces}. This terminology was motivated by the fact that if $w_{\mathsf{m}%
}=1 $ for every $\mathsf{m}$, $\Phi \left( t\right) =t^{2}$ and $k\in 
\mathbb{N}$, these spaces reduce to the set of all complex sequences\textsf{%
\ }$\left( p_{\mathsf{m}}\right) $ satisfying%
\begin{equation*}
\sum_{\mathsf{m}\in \mathbb{Z}}\left( 1+\mathsf{m}^{2}\right) ^{k}\left\vert
p_{\mathsf{m}}\right\vert ^{2}<\infty ,
\end{equation*}%
which thus emerge as the Fourier coefficients of periodic functions on the
real line belonging to a Sobolev space of order $k$. For further information
on this issue we refer the reader to \cite{Vuillermotbis}. As we shall see
in the following section, the results regarding the existence of continuous
and compact embeddings we prove there are natural generalizations of the
corresponding statements in Section 2 of \cite{Vuillermotter} relative to $%
\mathsf{s}_{k,\Phi ,w}$, $l_{k,\Phi ,w}$, $\mathsf{s}_{\Phi ,w}$, and $%
l_{\Phi ,w}$.

\section[Continuous and compact embeddings]{Continuous and compact embeddings for $\mathsf{s}_{\hat{\Phi}}$ and 
$l_{\hat{\Phi}}$}

Throughout this section we denote by $\Phi ^{-1}$, $\Psi ^{-1}$, $\Phi _{%
\mathsf{m}}^{-1}$, $\Psi _{\mathsf{m}}^{-1}$ the monotone inverses of the
Orlicz functions $\Phi $, $\Psi $, $\Phi _{\mathsf{m}}$, $\Psi _{\mathsf{m}}$%
, respectively. We begin with the following:

\bigskip

\textbf{Theorem 4. }\textit{Let }$\left( \Phi _{\mathsf{m}}\right) $\textit{%
\ and }$\left( \Psi _{\mathsf{m}}\right) $\textit{\ be two sequences of
Orlicz functions and let }$\left( \tau _{\mathsf{m}}\right) $\textit{\ be
the numerical sequence determined by the condition }$\Psi _{\mathsf{m}%
}\left( \tau _{\mathsf{m}}\right) =1$, \textit{which we assume to satisfy} 
\begin{equation}
\tau _{\mathsf{m}}\leq \tau <\infty  \label{boundedness}
\end{equation}%
\textit{for every }$\mathsf{m}$ \textit{and some }$\tau >0$. \textit{Then
the following statements hold:}

\textit{(a) If there exist constants }$\gamma ,t_{0}\in \left( 0,\infty
\right) $ \textit{such that the inequality }%
\begin{equation}
\Phi _{\mathsf{m}}\left( t\right) \leq \Psi _{\mathsf{m}}\left( \gamma
t\right)  \label{inequality}
\end{equation}%
\textit{holds for each }$t\in \left[ 0,t_{0}\right] $\textit{\ and every }$%
\mathsf{m}$, \textit{then there exists the continuous embedding}%
\begin{equation}
l_{\hat{\Psi}}\rightarrow l_{\hat{\Phi}}.  \label{embedding1}
\end{equation}%
\textit{(b) In addition to the above hypotheses, let us assume that}%
\begin{equation}
\liminf_{\left\vert \mathsf{m}\right\vert \rightarrow \infty }\Psi _{\mathsf{%
m}}\left( t\right) >0  \label{constraint}
\end{equation}%
\textit{for every }$t>0$. \textit{Then there exists the continuous embedding}
\begin{equation}
\mathsf{s}_{\hat{\Psi}}\rightarrow \mathsf{s}_{\hat{\Phi}}.
\label{embedding2}
\end{equation}

\bigskip

\textbf{Proof. }Let us prove (\ref{embedding1}). According to (\ref{norm}),
it is sufficient to show that there exists a finite constant $c_{\gamma
,t_{0,}\tau }>0$ depending solely on $\gamma ,t_{0}$ and $\tau $ such that
the inequality%
\begin{equation*}
\sum_{\mathsf{m}\in \mathbb{Z}}\Phi _{\mathsf{m}}\left( \frac{\left\vert p_{%
\mathsf{m}}\right\vert }{c_{\gamma ,t_{0},\tau }\left\Vert \mathsf{p}%
\right\Vert _{\hat{\Psi}}}\right) \leq 1
\end{equation*}%
holds for every non-zero $\mathsf{p}\in l_{\hat{\Psi}}$. Since we have%
\begin{equation*}
\sum_{\mathsf{m}\in \mathbb{Z}}\Psi _{\mathsf{m}}\left( \frac{\left\vert p_{%
\mathsf{m}}\right\vert }{\left\Vert \mathsf{p}\right\Vert _{\hat{\Psi}}}%
\right) \leq 1
\end{equation*}%
from (\ref{importantrelation}), Relation (\ref{boundedness}) then implies
the inequalities 
\begin{equation*}
\frac{\left\vert p_{\mathsf{m}}\right\vert }{c\left\Vert \mathsf{p}%
\right\Vert _{\hat{\Psi}}}\leq \frac{\Psi _{\mathsf{m}}^{-1}\left( 1\right) 
}{c}\leq \frac{\tau }{c}
\end{equation*}%
for each finite $c>0$ and every $\mathsf{m}\in \mathbb{Z}$, so that it is
sufficient to choose%
\begin{equation*}
c\geq \max \left( \gamma ,\frac{\tau }{t_{0}}\right) .
\end{equation*}%
Indeed, using (\ref{inequality}) we then have successively%
\begin{eqnarray*}
&&\sum_{\mathsf{m\in }\mathbb{Z}}\Phi _{\mathsf{m}}\left( \frac{\left\vert
p_{\mathsf{m}}\right\vert }{c\left\Vert \mathsf{p}\right\Vert _{\hat{\Psi}}}%
\right) \\
&\leq &\sum_{\mathsf{m\in }\mathbb{Z}}\Psi _{\mathsf{m}}\left( \frac{\gamma
\left\vert p_{\mathsf{m}}\right\vert }{c\left\Vert \mathsf{p}\right\Vert _{%
\hat{\Psi}}}\right) \leq \sum_{\mathsf{m\in }\mathbb{Z}}\Psi _{\mathsf{m}%
}\left( \frac{\left\vert p_{\mathsf{m}}\right\vert }{\left\Vert \mathsf{p}%
\right\Vert _{\hat{\Psi}}}\right) \leq 1,
\end{eqnarray*}%
as required.

In order to prove (\ref{embedding2}), it is sufficient to show that%
\begin{equation*}
\mathsf{s}_{\hat{\Psi}}\subseteq \mathsf{s}_{\hat{\Phi}}
\end{equation*}%
since the inclusion $\mathsf{s}_{\hat{\Psi}}\subseteq $\ $l_{\hat{\Psi}}$
and (\ref{embedding1}) hold. Thus, if $\mathsf{p}\in \mathsf{s}_{\hat{\Psi}}$
we have in particular 
\begin{equation}
\lim_{\left\vert \mathsf{m}\right\vert \rightarrow \infty }\Psi _{\mathsf{m}%
}\left( \frac{\left\vert p_{\mathsf{m}}\right\vert }{\rho }\right) =0
\label{limitbis}
\end{equation}%
for each $\rho \in \left( 0,\infty \right) $, so that%
\begin{equation*}
\lim_{\left\vert \mathsf{m}\right\vert \rightarrow \infty }\left( \frac{%
\left\vert p_{\mathsf{m}}\right\vert }{\rho }\right) =0.
\end{equation*}%
Indeed, if the preceding relation were not true there would exist $\epsilon
>0$ and a subsequence indexed by $\left( \mathsf{m}_{\mathsf{n}}\right) _{%
\mathsf{n}\in \mathbb{N}}$ such that%
\begin{equation*}
\frac{\left\vert p_{\mathsf{m}_{\mathsf{n}}}\right\vert }{\rho }\geq
\epsilon \text{,}
\end{equation*}%
thereby implying%
\begin{equation*}
0=\lim_{\mathsf{n}\rightarrow \infty }\Psi _{\mathsf{m}_{\mathsf{n}}}\left( 
\frac{\left\vert p_{\mathsf{m}_{\mathsf{n}}}\right\vert }{\rho }\right) \geq
\lim_{\mathsf{n}\rightarrow \infty }\Psi _{\mathsf{m}_{\mathsf{n}}}\left(
\epsilon \right) =0
\end{equation*}%
as a consequence of (\ref{limitbis}), thus contradicting (\ref{constraint}).
Therefore, there exists $\mathsf{M}_{t_{0}}\in \mathbb{N}$ such that $\rho
^{-1}\left\vert p_{\mathsf{m}}\right\vert \leq t_{0}$ whenever $\left\vert 
\mathsf{m}\right\vert \geq \mathsf{M}_{t_{0}}$ and so we get%
\begin{equation*}
\sum_{\left\vert \mathsf{m}\right\vert \geq \mathsf{M}_{t_{0}}}\Phi _{%
\mathsf{m}}\left( \frac{\left\vert p_{\mathsf{m}}\right\vert }{\rho }\right)
\leq \sum_{\left\vert \mathsf{m}\right\vert \geq \mathsf{M}_{t_{0}}}\Psi _{%
\mathsf{m}}\left( \frac{\left\vert p_{\mathsf{m}}\right\vert }{\gamma
^{-1}\rho }\right) <\infty
\end{equation*}%
for every $\rho \in \left( 0,\infty \right) $ from (\ref{inequality}), which
proves the desired inclusion. \ \ $\blacksquare $

\bigskip

Theorem 4 has the following consequence, showing that it provides a genuine
generalization to some of the results stated in Theorem 3 of \cite%
{Vuillermotter} as mentioned earlier:

\bigskip

\textbf{Corollary 1. }\textit{Let us consider the sequences }$\left( \Phi _{%
\mathsf{m}}\right) $\textit{\ and }$\left( \Psi _{\mathsf{m}}\right) $%
\textit{\ defined by}%
\begin{equation}
\Phi _{\mathsf{m}}(t):=w_{\mathsf{m}}\Phi (t)  \label{sequencefirst}
\end{equation}%
\textit{and }%
\begin{equation}
\Psi _{\mathsf{m}}(t):=w_{\mathsf{m}}\left( 1+\Psi \left( \left\vert \mathsf{%
m}\right\vert \right) \right) ^{k}\Psi (t)  \label{sequence}
\end{equation}%
\textit{respectively, where }$\Phi $\textit{, }$\Psi $\textit{\ stand for
fixed Orlicz functions and }$w=\left( w_{\mathsf{m}}\right) $\textit{, }$%
k\in \mathbb{R}$\textit{\ are as in (\ref{exempleparticulier}). Let us
assume that the inequality}%
\begin{equation}
\Phi \left( t\right) \leq \Psi \left( \gamma t\right)  \label{inequalityonze}
\end{equation}%
\textit{holds for each\ }$t\in \left[ 0,t_{0}\right] $\textit{\ for some }$%
\gamma ,t_{0}\in \left( 0,\infty \right) $\textit{. If \ \ }%
\begin{equation}
k\geq 0\text{, \ \ \ }\underline{w}:=\inf_{\mathsf{m}\in \mathbb{Z}}w_{%
\mathsf{m}}>0,  \label{conditionbis}
\end{equation}%
\textit{then there exist the continuous embeddings}%
\begin{equation*}
l_{k,\Psi ,w}\rightarrow l_{\Phi ,w}
\end{equation*}%
\textit{and}%
\begin{equation*}
\mathsf{s}_{k,\Psi ,w}\rightarrow \mathsf{s}_{\Phi ,w}.
\end{equation*}

\textbf{Proof. }It is plain that (\ref{sequencefirst}), (\ref{sequence}) and
(\ref{inequalityonze}) imply (\ref{inequality}) since $k\geq 0$, so that we
need only prove the boundedness of the sequence $\left( \tau _{\mathsf{m}%
}\right) $ and (\ref{constraint}). For convenience, let us define momentarily%
\begin{equation*}
\mu _{_{k,\Psi ,\mathsf{m}}}:=w_{\mathsf{m}}\left( 1+\Psi \left( \left\vert 
\mathsf{m}\right\vert \right) \right) ^{k}.
\end{equation*}%
The $\tau _{\mathsf{m}}$ are then determined by the condition $\tau _{%
\mathsf{m}}=\Psi ^{-1}\left( \mu _{_{k,\Psi ,\mathsf{m}}}^{-1}\right) $. But
(\ref{conditionbis}) implies%
\begin{equation*}
\mu _{_{k,\Psi ,\mathsf{m}}}\geq \underline{w},
\end{equation*}%
so that the boundedness%
\begin{equation*}
\tau _{\mathsf{m}}\leq \Psi ^{-1}\left( \underline{w}^{-1}\right) <\infty
\end{equation*}%
follows. On the other hand we have%
\begin{equation*}
\liminf_{\left\vert \mathsf{m}\right\vert \rightarrow \infty }\Psi _{\mathsf{%
m}}\left( t\right) \geq \underline{w}\Psi (t)>0
\end{equation*}%
for every $t>0$. \ \ $\blacksquare $

\bigskip

\textsc{Remarks.} (1) Let us consider the Orlicz functions%
\begin{eqnarray*}
\Phi _{\mathsf{m}}(t) &=&t^{\alpha _{\mathsf{m}}}, \\
\Psi _{\mathsf{m}}(t) &=&t^{\beta _{\mathsf{m}}}
\end{eqnarray*}%
with $\alpha _{\mathsf{m}},\beta _{\mathsf{m}}\in \left[ 1,\infty \right) $
and $\alpha _{\mathsf{m}}>\beta _{\mathsf{m}}$ for each $\mathsf{m}$. If we
set $t_{0}=\gamma ^{\frac{\beta _{\mathsf{m}}}{\alpha _{\mathsf{m}}-\beta _{%
\mathsf{m}}}}$ for any finite $\gamma \in \left( 0,\infty \right) $, then (%
\ref{inequality}) holds. Moreover, in this case we have $\tau _{\mathsf{m}%
}=\Psi _{\mathsf{m}}^{-1}(1)=1$ for every $\mathsf{m}$, a constant sequence.
Therefore Part (a) of Theorem 4 applies, so that we have the continuous
embedding%
\begin{equation*}
l_{\hat{\Psi}}\rightarrow l_{\hat{\Phi}}
\end{equation*}%
of the corresponding Nakano spaces. Moreover, if for instance $%
\lim_{\left\vert \mathsf{m}\right\vert \rightarrow \infty }\beta _{\mathsf{m}%
}=\beta $ then $\lim_{\left\vert \mathsf{m}\right\vert \rightarrow \infty
}t^{\beta _{\mathsf{m}}}=t^{\beta }>0$ for any $t>0$ so that%
\begin{equation*}
\mathsf{s}_{\hat{\Psi}}\rightarrow \mathsf{s}_{\hat{\Phi}}
\end{equation*}%
according to Part (b). For further extensions around this problem we refer
the reader to Section 3.3 in Chapter 3 of \cite{Diharu}.

(2) A natural question is whether the boundedness condition (\ref%
{boundedness}) and (\ref{constraint}) are also necessary for embeddings (\ref%
{embedding1}) and (\ref{embedding2}) to exist. The answer is clearly
negative, for if we consider for instance the functions%
\begin{equation*}
\Phi _{\mathsf{m}}\left( t\right) =\Psi _{\mathsf{m}}\left( t\right) =\frac{%
t^{\alpha }}{\left\vert \mathsf{m}\right\vert +1}
\end{equation*}%
where $\alpha \in \left[ 1,\infty \right) $, inequality (\ref{inequality})
is satisfied with $\gamma =1$, $t_{0}$ arbitrary, and the sequence $\left(
\tau _{\mathsf{m}}\right) $ is determined by $\tau _{\mathsf{m}}=\Psi _{%
\mathsf{m}}^{-1}\left( 1\right) =\left( \left\vert \mathsf{m}\right\vert
+1\right) ^{\frac{1}{\alpha }}$, which is clearly unbounded. Moreover,%
\begin{equation*}
\liminf_{\left\vert \mathsf{m}\right\vert \rightarrow \infty }\Psi _{\mathsf{%
m}}\left( t\right) =0
\end{equation*}%
for every $t>0$. Nonetheless, the indentity map in $l_{\hat{\Psi}}=\mathsf{s}%
_{\hat{\Psi}}=$ $\mathsf{s}_{\hat{\Phi}}=$ $l_{\hat{\Phi}}$ is continuous.

\bigskip

We proceed with the following result:

\bigskip

\textbf{Theorem 5. }\textit{Let us consider the sequence }$\left( \Phi _{%
\mathsf{m}}\right) $\textit{\ defined by }%
\begin{equation}
\Phi _{\mathsf{m}}(t):=\mu _{\left\vert \mathsf{m}\right\vert }\Phi \left(
t\right)  \label{definitionbis}
\end{equation}%
\textit{where }$\Phi $\textit{\ is a fixed Orlicz function and }$\left( \mu
_{\left\vert \mathsf{m}\right\vert }\right) $\textit{\ a numerical sequence
satisfying }$\mu _{\left\vert \mathsf{m}\right\vert }>0$\textit{\ for each }$%
\mathsf{m}$.\textit{\ Let }$\left( \Psi _{\mathsf{m}}\right) $ \textit{be a
sequence of} \textit{Orlicz functions such that the inequality}%
\begin{equation}
\Phi \left( t\right) \leq \Psi _{\mathsf{m}}\left( t\right)
\label{inequalityoitavo}
\end{equation}%
\textit{holds for each }$t\in \left( 0,\infty \right) $\textit{\ and every }$%
\mathsf{m}$\textit{. If the sequence }$\left( \mu _{\left\vert \mathsf{m}%
\right\vert }\right) $\textit{\ is bounded, then there exist the continuous
embeddings}%
\begin{equation}
l_{\hat{\Psi}}\rightarrow l_{\hat{\Phi}}  \label{embedding3}
\end{equation}%
\textit{and}%
\begin{equation}
\mathsf{s}_{\hat{\Psi}}\rightarrow \mathsf{s}_{\hat{\Phi}}.
\label{embedding6}
\end{equation}

\bigskip

\textbf{Proof. }By hypothesis there is a finite constant $c>0$ such that $%
\mu _{\left\vert \mathsf{m}\right\vert }\leq c$ for every $\mathsf{m}$. Then
for each $\mathsf{p}\in l_{\hat{\Psi}}$ there exists $\rho ^{\ast }>0$ such
that%
\begin{equation}
\sum_{\mathsf{m\in }\mathbb{Z}}\Phi _{\mathsf{m}}\left( \frac{\left\vert p_{%
\mathsf{m}}\right\vert }{\rho ^{\ast }}\right) \leq \sum_{\mathsf{m\in }%
\mathbb{Z}}\left( c\Phi \right) \left( \frac{\left\vert p_{\mathsf{m}%
}\right\vert }{\rho ^{\ast }}\right) \leq c\sum_{\mathsf{m\in }\mathbb{Z}%
}\Psi _{\mathsf{m}}\left( \frac{\left\vert p_{\mathsf{m}}\right\vert }{\rho
^{\ast }}\right) <\infty  \label{inequalitynono}
\end{equation}%
according to (\ref{definitionbis}) and (\ref{inequalityoitavo}), so that the
inclusion $l_{\hat{\Psi}}\subseteq l_{\hat{\Phi}}$ holds. Moreover, if $%
c\leq 1$ we have%
\begin{equation*}
\left\{ \rho >0:\sum_{\mathsf{m\in }\mathbb{Z}}\Psi _{\mathsf{m}}\left( 
\frac{\left\vert p_{\mathsf{m}}\right\vert }{\rho }\right) \leq 1\right\}
\subseteq \left\{ \rho >0:\sum_{\mathsf{m\in }\mathbb{Z}}\Phi _{\mathsf{m}%
}\left( \frac{\left\vert p_{\mathsf{m}}\right\vert }{\rho }\right) \leq
1\right\}
\end{equation*}%
since inequalities (\ref{inequalitynono}) also hold for those sufficiently
large $\rho >0$ that make the preceding sets non-empty (see the proof of
Theorem 1). Therefore we get%
\begin{equation}
\left\Vert \mathsf{p}\right\Vert _{\hat{\Phi}}\leq \left\Vert \mathsf{p}%
\right\Vert _{\hat{\Psi}}  \label{inequalitydecimo}
\end{equation}%
as a consequence of (\ref{norm}). Now if $c$ $>1$ we have%
\begin{equation*}
\Phi \left( \frac{1}{c}t\right) \leq \frac{1}{c}\Phi \left( t\right)
\end{equation*}%
for every $t\in \left( 0,\infty \right) $ from the convexity of $\Phi $
since $\Phi \left( 0\right) =0$, so that%
\begin{equation*}
\sum_{\mathsf{m\in }\mathbb{Z}}\left( c\Phi \right) \left( \frac{1}{c}\frac{%
\left\vert p_{\mathsf{m}}\right\vert }{\rho }\right) \leq \sum_{\mathsf{m\in 
}\mathbb{Z}}\Phi \left( \frac{\left\vert p_{\mathsf{m}}\right\vert }{\rho }%
\right) \leq \sum_{\mathsf{m\in }\mathbb{Z}}\Psi _{\mathsf{m}}\left( \frac{%
\left\vert p_{\mathsf{m}}\right\vert }{\rho }\right)
\end{equation*}%
as well for those $\rho >0$. Consequently, the inclusion%
\begin{equation}
\left\{ \rho >0:\sum_{\mathsf{m\in }\mathbb{Z}}\Psi _{\mathsf{m}}\left( 
\frac{\left\vert p_{\mathsf{m}}\right\vert }{\rho }\right) \leq 1\right\}
\subseteq \left\{ \rho >0:\sum_{\mathsf{m\in }\mathbb{Z}}\left( c\Phi
\right) \left( \frac{1}{c}\frac{\left\vert p_{\mathsf{m}}\right\vert }{\rho }%
\right) \leq 1\right\}  \label{inclusion}
\end{equation}%
also holds. Writing $\left\Vert .\right\Vert _{c\Phi }$ for the gauge norm (%
\ref{norm}) relative to the modular sequence space generated by $c\Phi $, we
then infer from (\ref{inclusion}) and the first inequality in (\ref%
{inequalitynono}) with $\rho $ instead of $\rho ^{\ast }$ that%
\begin{equation}
\left\Vert \mathsf{p}\right\Vert _{\hat{\Phi}}\leq \left\Vert \mathsf{p}%
\right\Vert _{c\Phi }\leq c\left\Vert \mathsf{p}\right\Vert _{\hat{\Psi}}.
\label{inequalitydecimoprimeiro}
\end{equation}%
Thus, defining%
\begin{equation*}
\hat{c}:=\max \left( 1,c\right)
\end{equation*}%
and taking (\ref{inequalitydecimo}) and (\ref{inequalitydecimoprimeiro})
into account we get 
\begin{equation*}
\left\Vert \mathsf{p}\right\Vert _{\hat{\Phi}}\leq \hat{c}\left\Vert \mathsf{%
p}\right\Vert _{\hat{\Psi}}
\end{equation*}%
for every $\mathsf{p}\in l_{\hat{\Psi}}$, so that (\ref{embedding3}) indeed
holds. The proof of (\ref{embedding6}) follows the same lines. \ \ $%
\blacksquare $

\bigskip

\textsc{Remark.} Intuitively, it is tempting to believe that the converse
statement also holds, namely, that if the sequence $\left( \mu _{\left\vert 
\mathsf{m}\right\vert }\right) $ is unbounded then there are no continuous
embeddings of the form (\ref{embedding3}) or (\ref{embedding6}). Whereas
this is true for a specific subclass of functions satisfying the conditions
of Theorem 5 (see Remark (2) below), such a statement is false in general.
This will be a consequence of the following result, which is of independent
interest:

\bigskip

\textbf{Proposition 2.}\textit{\ Let the sequence }$\left( \Phi _{\mathsf{m}%
}\right) $\textit{\ be of the form}%
\begin{equation*}
\Phi _{\mathsf{m}}(t)=\mu _{\left\vert \mathsf{m}\right\vert }t^{\alpha }
\end{equation*}%
\textit{for every }$\mathsf{m}$\textit{\ and }$\alpha \in \left[ 1,\infty
\right) $\textit{. Let }$\left( \mu _{\left\vert \mathsf{m}\right\vert
}\right) $ \textit{be\ any numerical sequence satisfying }$\mu _{\left\vert 
\mathsf{m}\right\vert }>0$\textit{\ and let }$\left( \Psi _{\mathsf{m}%
}\right) $ \textit{be\ a sequence of Orlicz functions. If }%
\begin{equation*}
\sum_{\mathsf{m\in }\mathbb{Z}}\mu _{\left\vert \mathsf{m}\right\vert }\left[
\Psi _{\mathsf{m}}^{-1}\left( 1\right) \right] ^{\alpha }<\infty ,
\end{equation*}%
\textit{then there exist the continuous embeddings}%
\begin{equation}
l_{\hat{\Psi}}\rightarrow l_{\hat{\Phi}}  \label{embedding7}
\end{equation}%
\textit{and}%
\begin{equation}
\mathsf{s}_{\hat{\Psi}}\rightarrow \mathsf{s}_{\hat{\Phi}}.
\label{embedding8}
\end{equation}

\bigskip

\textbf{Proof. }Since $\Psi _{\mathsf{m}}\left( \frac{\left\vert p_{\mathsf{m%
}}\right\vert }{\rho }\right) \geq 0$ for each $\mathsf{m}\in \mathbb{Z}$
and every finite $\rho >0$, it is sufficient to note that the inclusion%
\begin{equation*}
\left\{ \rho >0:\sum_{\mathsf{m}\in \mathbb{Z}}\Psi _{\mathsf{m}}\left( 
\frac{\left\vert p_{\mathsf{m}}\right\vert }{\rho }\right) \leq 1\right\}
\subseteq \left\{ \rho >0:\rho \geq \frac{\left\vert p_{\mathsf{m}%
}\right\vert }{\Psi _{\mathsf{m}}^{-1}\left( 1\right) }\right\}
\end{equation*}%
holds, from which we obtain%
\begin{equation*}
\left\vert p_{\mathsf{m}}\right\vert \leq \Psi _{\mathsf{m}}^{-1}\left(
1\right) \left\Vert \mathsf{p}\right\Vert _{\hat{\Psi}}
\end{equation*}%
\bigskip following (\ref{norm}). Therefore we get 
\begin{equation*}
\left\Vert \mathsf{p}\right\Vert _{\hat{\Phi}}^{\alpha }=\sum_{\mathsf{m}\in 
\mathbb{Z}}\mu _{\left\vert \mathsf{m}\right\vert }\left\vert p_{\mathsf{m}%
}\right\vert ^{\alpha }\leq \sum_{\mathsf{m}\in \mathbb{Z}}\mu _{\left\vert 
\mathsf{m}\right\vert }\left[ \Psi _{\mathsf{m}}^{-1}\left( 1\right) \right]
^{\alpha }\times \left\Vert \mathsf{p}\right\Vert _{\hat{\Psi}}^{\alpha },
\end{equation*}%
which is the desired result. \ \ $\blacksquare $

\bigskip

\textsc{Remarks.} (1) As a concrete example, let us consider the sequences $%
\left( \Phi _{\mathsf{m}}\right) $\textit{\ }and $\left( \Psi _{\mathsf{m}%
}\right) $ given by%
\begin{equation}
\Phi _{\mathsf{m}}(t)=\mu _{\left\vert \mathsf{m}\right\vert }t
\label{exempleparticulierbis}
\end{equation}%
and%
\begin{equation}
\Psi _{\mathsf{m}}\left( t\right) =\nu _{\left\vert \mathsf{m}\right\vert
}\left( \exp \left[ t\right] -1\right) ,  \label{exempleparticulierter}
\end{equation}%
respectively, where $\left( \mu _{\left\vert \mathsf{m}\right\vert }\right) $
is any numerical sequence satisfying $\mu _{\left\vert \mathsf{m}\right\vert
}>0$ and $\left( \nu _{\left\vert \mathsf{m}\right\vert }\right) $ a
numerical sequence satisfying $\nu _{\left\vert \mathsf{m}\right\vert }\geq
1 $ for every $\mathsf{m}$. If%
\begin{equation}
\sum_{\mathsf{m\in }\mathbb{Z}}\mu _{\left\vert \mathsf{m}\right\vert }\nu
_{\left\vert \mathsf{m}\right\vert }^{-1}<\infty ,  \label{conditionquinto}
\end{equation}%
then there exist the continuous embeddings (\ref{embedding7}) and (\ref%
{embedding8}).\textit{\ }Indeed we have%
\begin{equation*}
\Psi _{\mathsf{m}}^{-1}\left( 1\right) =\ln \left( 1+\nu _{\left\vert 
\mathsf{m}\right\vert }^{-1}\right)
\end{equation*}%
and thereby%
\begin{equation*}
\sum_{\mathsf{m\in }\mathbb{Z}}\mu _{\left\vert \mathsf{m}\right\vert }\ln
\left( 1+\nu _{\left\vert \mathsf{m}\right\vert }^{-1}\right) \leq \sum_{%
\mathsf{m\in }\mathbb{Z}}\mu _{\left\vert \mathsf{m}\right\vert }\nu
_{\left\vert \mathsf{m}\right\vert }^{-1}<\infty ,
\end{equation*}%
so that those embeddings follow from Proposition 2. But we also note that (%
\ref{exempleparticulierbis}) and (\ref{exempleparticulierter}) imply (\ref%
{inequalityoitavo}). Therefore, Condition (\ref{conditionquinto}) shows that
we can compensate for the possible unboundedness of $\left( \mu _{\left\vert 
\mathsf{m}\right\vert }\right) $ by choosing a $\left( \nu _{\left\vert 
\mathsf{m}\right\vert }\right) $ that grows sufficiently fast in order to
get continuous embeddings of the above form. Of course examples abound, for
instance 
\begin{equation*}
\Phi _{\mathsf{m}}(t)=\left( \mathsf{m}^{2}+1\right) t
\end{equation*}%
and%
\begin{equation*}
\Psi _{\mathsf{m}}\left( t\right) =\left( \mathsf{m}^{4}+1\right) \left(
\exp \left[ t\right] -1\right) .
\end{equation*}%
Consequently, the converse statement to Theorem 5 does not hold in general.

(2) The preceding remark notwithstanding, there is at least a small subclass
of modular spaces for which the converse statement of Theorem 5 is true.
Indeed, if the sequence $\left( \mu _{\left\vert \mathsf{m}\right\vert
}\right) $ is unbounded there exists a subsequence $\left( \mu _{\left\vert 
\mathsf{m}_{\mathsf{n}}\right\vert }\right) $ such that $\mu _{\left\vert 
\mathsf{m}_{\mathsf{n}}\right\vert }>\mathsf{n}$ for every $\mathsf{n}\in 
\mathbb{N}$. Then, let $\left( \Phi _{\mathsf{m}}\right) $ and $\left( \Psi
_{\mathsf{m}}\right) $ be of the form%
\begin{eqnarray*}
\Phi _{\mathsf{m}}(t) &=&\mu _{\left\vert \mathsf{m}\right\vert }\Phi (t), \\
\Psi _{\mathsf{m}}\left( t\right) &=&\Phi (t)
\end{eqnarray*}%
for every $\mathsf{m}$ where $\Phi $ is an Orlicz function, so that (\ref%
{inequalityoitavo}) holds. Let us furthermore assume that $\Phi $ satisfies
the global $\Delta _{2}$-condition in the form (\ref{inequalitysextobis}),
and choose $\mathsf{p=}\left( p_{\mathsf{m}}\right) $ as%
\begin{equation*}
p_{\mathsf{m}}=\left\{ 
\begin{array}{c}
\text{ }\Phi ^{-1}\left( \frac{1}{\mathsf{n}^{2}}\right) \text{\ if \ }%
\mathsf{m}=\mathsf{m}_{\mathsf{n}}, \\ 
\\ 
0\text{ \ \ if \ }\mathsf{m}\neq \mathsf{m}_{\mathsf{n}}.%
\end{array}%
\right.
\end{equation*}%
Then we get 
\begin{equation*}
\sum_{\mathsf{m\in }\mathbb{Z}}\Psi _{\mathsf{m}}\left( \frac{\left\vert p_{%
\mathsf{m}}\right\vert }{\rho }\right) =\sum_{\mathsf{m=m}_{\mathsf{n}}}\Phi
\left( \frac{\left\vert p_{\mathsf{m}}\right\vert }{\rho }\right) \leq \eta
_{\rho }\sum_{\mathsf{n=1}}^{\infty }\frac{1}{\mathsf{n}^{2}}<\infty
\end{equation*}%
for every $\rho >0$ and some finite constant $\eta _{\rho }>0$. This means
that $\mathsf{p\in s}_{\hat{\Psi}}$, whereas%
\begin{equation*}
\sum_{\mathsf{m\in }\mathbb{Z}}\Phi _{\mathsf{m}}\left( \left\vert p_{%
\mathsf{m}}\right\vert \right) =\sum_{\mathsf{m=m}_{\mathsf{n}}}^{\infty }%
\frac{\mu _{\left\vert \mathsf{m}_{\mathsf{n}}\right\vert }}{\mathsf{n}^{2}}%
\geq \sum_{\mathsf{n=1}}^{\infty }\frac{1}{\mathsf{n}}=\infty ,
\end{equation*}%
so that $\mathsf{p\notin s}_{\hat{\Phi}}$. Therefore $\mathsf{s}_{\hat{\Psi}%
}\nsubseteq \mathsf{s}_{\hat{\Phi}}$, and \textit{a fortiori} there are no
continuous embeddings of the form (\ref{embedding7}) and (\ref{embedding8})
in this case since $l_{\hat{\Psi}}=\mathsf{s}_{\hat{\Psi}}\nsubseteq \mathsf{%
s}_{\hat{\Phi}}=l_{\hat{\Phi}}$.

\bigskip

We now turn to the investigation of the existence problem of compact
embeddings. Its solution is not as simple and direct as that regarding the
existence of continuous embeddings, the latter merely requiring pointwise
comparisons between the sequences $\left( \Phi _{\mathsf{m}}\right) $\textit{%
\ }and\textit{\ }$\left( \Psi _{\mathsf{m}}\right) $ as in (\ref{inequality}%
) or (\ref{inequalityoitavo}). We begin with the following characterization
of bounded sets $\mathcal{K}$ in $l_{\hat{\Phi}}$ or $\mathsf{s}_{\hat{\Phi}%
} $ in terms of the $\Phi _{\mathsf{m}}$:

\bigskip

\textbf{Lemma 1. }\textit{Let }$\kappa \in \left( 0,\infty \right) .$ 
\textit{Then we have }$\left\Vert \mathsf{p}\right\Vert _{\hat{\Phi}}\leq
\kappa $\textit{\ for every }$\mathsf{p}\in \mathcal{K}$ \textit{if, and
only if,}%
\begin{equation}
\sum_{\mathsf{m}\in \mathbb{Z}}\Phi _{\mathsf{m}}\left( \frac{\left\vert p_{%
\mathsf{m}}\right\vert }{\kappa }\right) \leq 1  \label{relationbis}
\end{equation}%
\textit{for such} $\mathsf{p}$.

\bigskip

\textbf{Proof.} The statement is obvious if $\mathsf{p=0}$, so let us assume
that $\mathsf{p\neq 0}$. If (\ref{relationbis}) holds for $\mathsf{p}\in 
\mathcal{K}$, then we have $\left\Vert \mathsf{p}\right\Vert _{\hat{\Phi}%
}\leq \kappa $ as a consequence of (\ref{norm}).\ Conversely, if $\left\Vert 
\mathsf{p}\right\Vert _{\hat{\Phi}}\leq \kappa $ we get%
\begin{equation*}
\sum_{\mathsf{m}\in \mathbb{Z}}\Phi _{\mathsf{m}}\left( \frac{\left\vert p_{%
\mathsf{m}}\right\vert }{\kappa }\right) \leq \sum_{\mathsf{m}\in \mathbb{Z}%
}\Phi _{\mathsf{m}}\left( \frac{\left\vert p_{\mathsf{m}}\right\vert }{%
\left\Vert \mathsf{p}\right\Vert _{\hat{\Phi}}}\right) \leq 1
\end{equation*}%
from (\ref{importantrelation}) since $\Phi _{\mathsf{m}}$ is monotone
increasing. \ \ $\blacksquare $

\bigskip

Our next preparatory result characterizes the relative compactness of
bounded sets in $\mathsf{s}_{\hat{\Phi}}$. Its proof shows that the given
arguments also hold for any bounded subset in $l_{\hat{\Phi}}$ as the
specific structure of $\mathsf{s}_{\hat{\Phi}}$ and $l_{\hat{\Phi}}$ is
irrelevant there (in fact, the argument holds in any Banach space, see also 
\cite{Bisgard}). However, we shall omit further reference to $l_{\hat{\Phi}}$
in the sequel as the proof of our main compactness results will exemplify
the importance of the existence of a Schauder basis in the space under
consideration.

\bigskip

\textbf{Lemma 2. }\textit{Let }$\left( \Phi _{\mathsf{m}}\right) $\textit{\
be a sequence of Orlicz functions and let }$\mathsf{s}_{\hat{\Phi}}$\textit{%
\ be the corresponding small modular sequence space. Then a bounded set }$%
\mathcal{K}\subset \mathsf{s}_{\hat{\Phi}}$\textit{\ is relatively compact
if, and only if, for each }$\epsilon \in \left( 0,\infty \right) $\textit{\
there exists a finite-dimensional subspace }$\mathcal{S}_{\epsilon }\subset 
\mathsf{s}_{\hat{\Phi}}$\textit{\ such that }%
\begin{equation}
\operatorname{dist}\left( \mathsf{p},\mathcal{S}_{\epsilon }\right) <\epsilon
\label{distance}
\end{equation}%
\textit{for every }$\mathsf{p}\in \mathcal{\bar{K}}$\textit{, the closure of 
}$\mathcal{K}$ \textit{in} $\mathsf{s}_{\hat{\Phi}}$.

\bigskip

\textbf{Proof. }Let us first assume that $\mathcal{K}$ is relatively compact
in $\mathsf{s}_{\hat{\Phi}}$. Then its closure $\mathcal{\bar{K}}$ is
totally bounded there, that is, for each $\epsilon \in \left( 0,\infty
\right) $ there exist finitely many open balls $\mathcal{B}\left( \epsilon ,%
\mathsf{c}_{j}\right) $ of radius $\epsilon $ and centers $\mathsf{c}_{j}\in 
\mathcal{\bar{K}}$ such that 
\begin{equation}
\mathcal{\bar{K}\subseteq \cup }_{j=1}^{N_{\epsilon }}\mathcal{B}\left(
\epsilon ,\mathsf{c}_{j}\right) .  \label{finitecovering}
\end{equation}%
It is therefore clear that (\ref{distance}) is satisfied if we choose%
\begin{equation*}
\mathcal{S}_{\epsilon }=\operatorname{span}\left\{ \mathsf{c}_{1},...,\mathsf{c}%
_{N_{\epsilon }}\right\} \text{.}
\end{equation*}%
Indeed, let us assume there exists a $\mathsf{p}^{\ast }\in \mathcal{\bar{K}}
$ such that (\ref{distance}) does not hold. Then we have \textit{a fortiori }%
\begin{equation*}
\left\Vert \mathsf{p}^{\ast }-\mathsf{c}_{j}\right\Vert _{\hat{\Phi}}\geq
\epsilon
\end{equation*}%
for every $j\in \left\{ 1,...,N_{\epsilon }\right\} $, thereby contradicting
(\ref{finitecovering}).

Conversely, let us assume that there exists a finite-dimensional subspace $%
\mathcal{S}_{\frac{\epsilon }{3}}\subset \mathsf{s}_{\hat{\Phi}}$ satisfying%
\begin{equation*}
\operatorname{dist}\left( \mathsf{p},\mathcal{S}_{\frac{\epsilon }{3}}\right) <\frac{%
\epsilon }{3}
\end{equation*}%
for each $\mathsf{p}\in \mathcal{\bar{K}}$. Then for any such $\mathsf{p}$
there exists a $\mathsf{q}^{\ast }\in \mathcal{S}_{\frac{\epsilon }{3}}$
such that 
\begin{equation}
\left\Vert \mathsf{p-q}^{\ast }\right\Vert _{\hat{\Phi}}<\frac{\epsilon }{3}.
\label{inequalitybis}
\end{equation}%
Let us now define the set%
\begin{equation*}
\mathcal{L}_{\epsilon }:=\left\{ \mathsf{q}\in \mathcal{S}_{\frac{\epsilon }{%
3}}:\operatorname{dist}\left( \mathsf{q},\mathcal{\bar{K}}\right) <\frac{\epsilon }{3%
}\right\} .
\end{equation*}%
Then for any $\mathsf{q}\in \mathcal{L}_{\epsilon }$ there exists a $\mathsf{%
p}^{\ast }\in \mathcal{\bar{K}}$ such that 
\begin{equation}
\left\Vert \mathsf{q}-\mathsf{p}^{\ast }\right\Vert _{\hat{\Phi}}<\frac{%
\epsilon }{3},  \label{inequalityter}
\end{equation}%
which entails the boundedness of $\mathcal{L}_{\epsilon }$ since $\mathcal{%
\bar{K}}$ is bounded. Indeed, since there exists $\kappa \in \left( 0,\infty
\right) $ such that $\left\Vert \mathsf{p}^{\ast }\right\Vert _{_{\hat{\Phi}%
}}\leq \kappa $ we have%
\begin{equation*}
\left\Vert \mathsf{q}\right\Vert _{\hat{\Phi}}\leq \left\Vert \mathsf{q}-%
\mathsf{p}^{\ast }\right\Vert _{\hat{\Phi}}+\left\Vert \mathsf{p}^{\ast
}\right\Vert _{_{\hat{\Phi}}}\leq \frac{\epsilon }{3}+\kappa .
\end{equation*}%
Therefore, $\mathcal{L}_{\epsilon }$ is relatively compact as a bounded set
in a finite-dimensional subspace of $\mathsf{s}_{\hat{\Phi}}$, so that there
exist finitely many open balls $\mathcal{B}\left( \frac{\epsilon }{3},%
\mathsf{\hat{c}}_{j}\right) $ of radius $\frac{\epsilon }{3}$ and centers $%
\mathsf{\hat{c}}_{j}\in \mathcal{\bar{L}}_{\epsilon }$ such that%
\begin{equation}
\mathcal{\bar{L}}_{\epsilon }\mathcal{\subseteq \cup }_{j=1}^{\hat{N}%
_{\epsilon }}\mathcal{B}\left( \frac{\epsilon }{3},\mathsf{\hat{c}}%
_{j}\right) .  \label{finitecoveringbis}
\end{equation}%
We complete the proof by showing that (\ref{finitecoveringbis}) allows one
to construct a finite covering of $\mathcal{\bar{K}}$. On the one hand, it
follows from (\ref{inequalityter}) that for every $j\in \left\{ 1,...,\hat{N}%
_{\epsilon }\right\} $ there exists $c_{j}\in \mathcal{\bar{K}}$ such that%
\begin{equation}
\left\Vert \mathsf{\hat{c}}_{j}-\mathsf{c}_{j}\right\Vert _{\hat{\Phi}}<%
\frac{\epsilon }{3}.  \label{inequalityquarto}
\end{equation}%
On the other hand, it is plain that the $\mathsf{q}^{\ast }$ that appears in
(\ref{inequalitybis}) belongs to $\mathcal{L}_{\epsilon }$. Thus, as a
consequence of (\ref{finitecoveringbis}) there exists $j^{\ast }\in \left\{
1,...,\hat{N}_{\epsilon }\right\} $ such that%
\begin{equation}
\mathsf{q}^{\ast }\in \mathcal{B}\left( \frac{\epsilon }{3},\mathsf{\hat{c}}%
_{j^{\ast }}\right) .  \label{oneball}
\end{equation}%
Then, for every $\mathsf{p}\in \mathcal{\bar{K}}$ we have%
\begin{eqnarray*}
&&\left\Vert \mathsf{p}-\mathsf{c}_{j^{\ast }}\right\Vert _{\hat{\Phi}} \\
&\leq &\left\Vert \mathsf{p-q}^{\ast }\right\Vert _{\hat{\Phi}}+\left\Vert 
\mathsf{q}^{\ast }-\mathsf{\hat{c}}_{j^{\ast }}\right\Vert _{\hat{\Phi}%
}+\left\Vert \mathsf{\hat{c}}_{j^{\ast }}-\mathsf{c}_{j^{\ast }}\right\Vert
_{\hat{\Phi}}<\epsilon
\end{eqnarray*}%
by virtue of (\ref{inequalitybis}), (\ref{inequalityquarto}) and (\ref%
{oneball}). Therefore we have%
\begin{equation*}
\mathsf{p}\in \mathcal{B}\left( \epsilon ,\mathsf{c}_{j^{\ast }}\right)
\subset \mathcal{\cup }_{j=1}^{\hat{N}_{\epsilon }}\mathcal{B}\left(
\epsilon ,\mathsf{c}_{j}\right) ,
\end{equation*}%
which proves that $\mathcal{K}$ is relatively compact in $\mathsf{s}_{\hat{%
\Phi}}$. \ \ $\blacksquare $

\bigskip

The proof of our first compactness result relies on the uniform global $%
\Delta _{2}$-condition at the origin and on the existence of a Schauder
basis in $\mathsf{s}_{\hat{\Phi}}=l_{\hat{\Phi}}$:

\bigskip

\textbf{Theorem 6. }\textit{Let }$\left( \Phi _{\mathsf{m}}\right) $\textit{%
\ and }$\left( \Psi _{\mathsf{m}}\right) $\textit{\ be two sequences of
Orlicz functions, and let us assume that }$\left( \Phi _{\mathsf{m}}\right) $
\textit{satisfies} \textit{the uniform global }$\Delta _{2}$\textit{%
-condition. If the condition }%
\begin{equation}
\lim_{\left\vert \mathsf{m}\right\vert \rightarrow \infty
}\sup_{0<t\leqslant \Psi _{\mathsf{m}}^{-1}\left( 1\right) }\frac{\Phi _{%
\mathsf{m}}\left( t\right) }{\Psi _{\mathsf{m}}\left( t\right) }=0
\label{condition}
\end{equation}%
\textit{holds, then there exists the compact embedding}%
\begin{equation}
\mathsf{s}_{\hat{\Psi}}\hookrightarrow \mathsf{s}_{\hat{\Phi}}.
\label{compactembeddingbis}
\end{equation}

\bigskip

\textbf{Proof. }Let $\mathcal{K}\subset \mathsf{s}_{\hat{\Psi}}$ be a
bounded set with $\left\Vert \mathsf{p}\right\Vert _{\hat{\Psi}}\leq \kappa
<\infty $ for each\textit{\ }$\mathsf{p}\in \mathcal{K}$. From Lemma 1 we
have%
\begin{equation}
\sum_{\mathsf{m}\in \mathbb{Z}}\Psi _{\mathsf{m}}\left( \frac{\left\vert p_{%
\mathsf{m}}\right\vert }{\kappa }\right) \leq 1  \label{relationter}
\end{equation}%
and \textit{a fortiori}%
\begin{equation}
\kappa ^{-1}\left\vert p_{\mathsf{m}}\right\vert \leq \Psi _{\mathsf{m}%
}^{-1}\left( 1\right)  \label{conditionsexto}
\end{equation}%
for every $\mathsf{m}\in \mathbb{Z}$. Furthermore, for any $\epsilon \in
\left( 0,\infty \right) $ we have%
\begin{equation}
\Phi _{\mathsf{m}}\left( \frac{2\left\vert p_{\mathsf{m}}\right\vert }{%
\epsilon }\right) =\Phi _{\mathsf{m}}\left( \frac{2\kappa }{\epsilon }\frac{%
\left\vert p_{\mathsf{m}}\right\vert }{\kappa }\right) \leq \eta _{\epsilon
,\kappa }\Phi _{\mathsf{m}}\left( \frac{\left\vert p_{\mathsf{m}}\right\vert 
}{\kappa }\right)  \label{inequalitydoze}
\end{equation}%
for some finite $\eta _{\epsilon ,\kappa }>0$ whenever $\left\vert \mathsf{m}%
\right\vert \geq \mathsf{M}$, as a consequence of the uniform global $\Delta
_{2}$-condition at the origin in the form (\ref{inequalitysextobis}) with $%
\theta =2\kappa \epsilon ^{-1}$. Moreover, from (\ref{condition}) and (\ref%
{conditionsexto}) we may write%
\begin{equation}
\Phi _{\mathsf{m}}\left( \frac{\left\vert p_{\mathsf{m}}\right\vert }{\kappa 
}\right) \leq \eta _{\epsilon ,\kappa }^{-1}\Psi _{\mathsf{m}}\left( \frac{%
\left\vert p_{\mathsf{m}}\right\vert }{\kappa }\right)
\label{inequalitytreize}
\end{equation}%
uniformly in $\mathsf{p}$ whenever $\left\vert \mathsf{m}\right\vert \geq 
\mathsf{M}_{\epsilon ,\kappa }^{\prime }$ for some $\mathsf{M}_{\epsilon
,\kappa }^{\prime }\in \mathbb{N}$. Defining then $\mathsf{M}_{\epsilon
,\kappa }:=\allowbreak\max \left( \mathsf{M},\mathsf{M}_{\epsilon ,\kappa }^{\prime
}\right) $ and combining (\ref{inequalitydoze}) with (\ref{inequalitytreize}%
), we obtain%
\begin{equation*}
\displaystyle\sum\limits_{\left\vert \mathsf{m}\right\vert \geq \mathsf{M}_{\epsilon
,\kappa }}\Phi _{\mathsf{m}}\left( \frac{2\left\vert p_{\mathsf{m}%
}\right\vert }{\epsilon }\right) \leq \displaystyle\sum\limits_{\left\vert \mathsf{m}%
\right\vert \geq \mathsf{M}_{\epsilon ,\kappa }}\Psi _{\mathsf{m}}\left( 
\frac{\left\vert p_{\mathsf{m}}\right\vert }{\kappa }\right) \leq 1
\end{equation*}%
from (\ref{relationter}) or, equivalently,

\begin{equation}
\frac{\epsilon }{2}\in \left\{ \rho >0:\sum_{\left\vert \mathsf{m}%
\right\vert \geq \mathsf{M}_{\epsilon ,\kappa }}\Phi _{\mathsf{m}}\left( 
\frac{\left\vert p_{\mathsf{m}}\right\vert }{\rho }\right) \leq 1\right\} .
\label{set}
\end{equation}%
Now since $\mathsf{s}_{\hat{\Psi}}\subseteq \mathsf{s}_{\hat{\Phi}}$ as a
consequence of (\ref{condition}), every $\mathsf{p}\in \mathcal{K}$ may be
viewed as an element of $\mathsf{s}_{\hat{\Phi}}$ so that it can be expanded
as%
\begin{equation*}
\mathsf{p=}\sum_{\mathsf{m}\in \mathbb{Z}}p_{\mathsf{m}}\mathsf{e}_{\mathsf{m%
}}
\end{equation*}%
in the norm topology of that space, where $\left( \mathsf{e}_{\mathsf{m}%
}\right) $ stands for the Schauder basis introduced in Theorem 2. Defining
now%
\begin{equation*}
\mathsf{p}_{\mathsf{M}_{\epsilon ,\kappa }}:\mathsf{=}\sum_{\left\vert 
\mathsf{m}\right\vert <\mathsf{M}_{\epsilon ,\kappa }}p_{\mathsf{m}}\mathsf{e%
}_{\mathsf{m}}
\end{equation*}%
we claim that 
\begin{equation}
\left\Vert \mathsf{p-p}_{\mathsf{M}_{\epsilon ,\kappa }}\right\Vert _{\hat{%
\Phi}}\leq \frac{\epsilon }{2}.  \label{estimatebis}
\end{equation}%
In order to see this it is sufficient to show that%
\begin{eqnarray*}
&&\left\{ \rho >0:\sum_{\left\vert \mathsf{m}\right\vert \geq \mathsf{M}%
_{\epsilon ,\kappa }}\Phi _{\mathsf{m}}\left( \frac{\left\vert p_{\mathsf{m}%
}\right\vert }{\rho }\right) \leq 1\right\} \\
&=&\left\{ \rho >0:\sum_{\mathsf{m\in }\mathbb{Z}}\Phi _{\mathsf{m}}\left( 
\frac{\left\vert \left( \mathsf{p-p}_{\mathsf{M}_{\epsilon ,\kappa }}\right)
_{\mathsf{m}}\right\vert }{\rho }\right) \leq 1\right\}
\end{eqnarray*}%
according to (\ref{norm}) and (\ref{set}). But this is an immediate
consequence of the fact that the $\mathsf{m}^{\text{th }}$component of $%
\mathsf{p-p}_{\mathsf{M}_{\epsilon ,\kappa }}$ reads 
\begin{equation*}
\left( \mathsf{p-p}_{\mathsf{M}_{\epsilon }}\right) _{\mathsf{m}}=\left\{ 
\begin{array}{c}
p_{\mathsf{m}}\text{ \ if \ }\mathsf{m}\in \mathbb{Z}\cap \left( -\infty ,-%
\mathsf{M}_{\epsilon ,\kappa }\right] \text{ or if \ }\mathsf{m}\in \mathbb{Z%
}\cap \left[ \mathsf{M}_{\epsilon ,\kappa },\infty \right) \text{\ } \\ 
\\ 
0\text{ \ otherwise.}%
\end{array}%
\right.
\end{equation*}%
Lastly, let us define the finite-dimensional subspace $\mathcal{S}_{\epsilon
,\kappa }\subset \mathsf{s}_{\hat{\Phi}}$ by%
\begin{equation*}
\mathcal{S}_{\epsilon ,\kappa }:=\operatorname{span}\left\{ \mathsf{e}_{_{-%
\mathsf{M}_{\epsilon ,\kappa }}}\text{,...., }\mathsf{e}_{_{\mathsf{M}%
_{\epsilon ,\kappa }}}\right\} .
\end{equation*}%
Then we clearly have $\mathsf{p}_{\mathsf{M}_{\epsilon ,\kappa }}\in 
\mathcal{S}_{\epsilon ,\kappa }$ and so for every $\mathsf{p}\in \mathcal{K}$
we get%
\begin{equation*}
\operatorname{dist}\left( \mathsf{p,}\mathcal{S}_{\epsilon ,\kappa }\right) \leq 
\frac{\epsilon }{2}
\end{equation*}%
according to (\ref{estimatebis}). It remains to obtain a similar estimate
for any $\mathsf{p}\in \mathcal{\bar{K}\diagdown K}$, in which case there
exists $\mathsf{p}_{\epsilon }\in \mathcal{K}$ such that $\left\Vert \mathsf{%
p-p}_{\epsilon }\right\Vert _{\hat{\Phi}}<$ $\frac{\epsilon }{2}$. Therefore
we have%
\begin{equation}
\operatorname{dist}\left( \mathsf{p,}\mathcal{S}_{\epsilon ,\kappa }\right) \leq 
\frac{\epsilon }{2}+\inf_{\mathsf{q}\in \mathcal{S}_{\epsilon ,\kappa
}}\left\Vert \mathsf{p}_{\epsilon }-\mathsf{q}\right\Vert _{\hat{\Phi}}\leq
\epsilon  \label{distancebis}
\end{equation}%
as desired, which proves the relative compactness of $\mathcal{K}$ according
to Lemma 2. \ \ $\blacksquare $

\bigskip

\textsc{Remark. }Whereas (\ref{condition}) is a sufficient condition to
ensure relative compactness, it is by no means necessary. Indeed, let $h_{%
\mathbb{C}}^{1}$ be the space of Sobolev type consisting of all complex
sequences $\mathsf{p=}\left( p_{\mathsf{m}}\right) $ satisfying%
\begin{equation*}
\sum_{\mathsf{m}\in \mathbb{Z}}\left( 1+\mathsf{m}^{2}\right) \left\vert p_{%
\mathsf{m}}\right\vert ^{2}<\infty ,
\end{equation*}%
and let $l_{\mathbb{C}}^{1}$ be the space of all complex, absolutely
convergent sequences. It is known from Proposition A.1 in Appendix A of \cite%
{Vuillermot} that the embedding%
\begin{equation}
h_{\mathbb{C}}^{1}\hookrightarrow l_{\mathbb{C}}^{1}
\label{compactembeddingter}
\end{equation}%
is compact. But $h_{\mathbb{C}}^{1}$ may be considered as the modular space
generated by the sequence of Orlicz functions%
\begin{equation*}
\Psi _{\mathsf{m}}\left( t\right) =\left( 1+\mathsf{m}^{2}\right) t^{2},
\end{equation*}%
while $l_{\mathbb{C}}^{1}$ corresponds to%
\begin{equation*}
\Phi _{\mathsf{m}}\left( t\right) =t.
\end{equation*}%
In this case it is clear that (\ref{condition}) does not hold.

\bigskip

The preceding result also leads to a characterization of the relative
compactness of bounded sets in $\mathsf{s}_{\hat{\Phi}}=l_{\hat{\Phi}}$,
particularly of those bounded sets in Nakano's spaces, thereby providing a
generalization of a standard compactness criterion for bounded sets in the
spaces $l_{\mathbb{C}}^{\alpha }$ consisting of all complex $\alpha $%
-summable sequences, where $\alpha \in \left[ 1,\infty \right) $:

\bigskip

\textbf{Corollary 2. }\textit{Let }$\left( \Phi _{\mathsf{m}}\right) $%
\textit{\ be a sequence of Orlicz functions that satisfy the uniform global }%
$\Delta _{2}$\textit{-condition, and let }$\mathcal{K}\subset \mathsf{s}_{%
\hat{\Phi}}$ \textit{be bounded. Then }$\mathcal{K}$ \textit{is relatively
compact if, and only if, the condition}%
\begin{equation}
\lim_{\mathsf{M}\rightarrow \infty }\sup_{\mathsf{p}\in \mathcal{K}%
}\sum_{\left\vert \mathsf{m}\right\vert \geq \mathsf{M}}\Phi _{\mathsf{m}%
}\left( \left\vert p_{\mathsf{m}}\right\vert \right) =0
\label{conditionquarto}
\end{equation}%
\textit{holds.}

\bigskip

\textbf{Proof. }It is plain that (\ref{condition}) implies (\ref%
{conditionquarto}). Arguing then as in the first part of the proof of
Theorem 6 by means of (\ref{inequalitydoze}) we infer the existence of a
constant $\mathsf{M}_{\epsilon ,\kappa }\in \mathbb{N}$ such that 
\begin{equation*}
\displaystyle\sum\limits_{\left\vert \mathsf{m}\right\vert \geq \mathsf{M}_{\epsilon
,\kappa }}\Phi _{\mathsf{m}}\left( \frac{2\left\vert p_{\mathsf{m}%
}\right\vert }{\epsilon }\right) \leq 1,
\end{equation*}%
thereby proving the relative compactness of $\mathcal{K}$.

Conversely, let $\mathcal{K}$ be relatively compact in $\mathsf{s}_{\hat{\Phi%
}}$ and let $\epsilon \in \left( 0,\infty \right) $. Then there exist
finitely many open balls $\mathcal{B}\left( \epsilon ,\mathsf{c}_{j}\right) $
of radius $\epsilon $ and centers $\mathsf{c}_{j}\in \mathcal{\bar{K}}$ such
that 
\begin{equation}
\mathcal{\bar{K}\subseteq \cup }_{j=1}^{N_{\epsilon }}\mathcal{B}\left(
\epsilon ,\mathsf{c}_{j}\right) .  \label{covering}
\end{equation}%
Moreover, since $\mathsf{c}_{j}\in \mathsf{s}_{\hat{\Phi}}$ for every $j\in
\left\{ 1,...,N_{\epsilon }\right\} $, there exists $\mathsf{M}_{\epsilon
,j}\in \mathbb{N}$ such that 
\begin{equation*}
\sum_{\left\vert \mathsf{m}\right\vert \geq \mathsf{M}_{\epsilon ,j}}\Phi _{%
\mathsf{m}}\left( \left\vert \mathsf{c}_{j,\mathsf{m}}\right\vert \right)
\leq \epsilon
\end{equation*}%
where $\mathsf{c}_{j,\mathsf{m}}$ stands for the $\mathsf{m}^{\mathsf{th}}$
component of $\mathsf{c}_{j}$. Defining%
\begin{equation}
\mathsf{M}_{\epsilon }:=\max_{j\in \left\{ 1,...,N_{\epsilon }\right\} }%
\mathsf{M}_{\epsilon ,j}  \label{maximum}
\end{equation}%
we then have \textit{a fortiori}%
\begin{equation}
\sum_{\left\vert \mathsf{m}\right\vert \geq \mathsf{M}_{\epsilon }}\Phi _{%
\mathsf{m}}\left( \left\vert \mathsf{c}_{j,\mathsf{m}}\right\vert \right)
\leq \epsilon  \label{smalltail}
\end{equation}%
uniformly in $j$. Let us now restrict ourselves momentarily to the case $%
\epsilon \in \left( 0,1\right] $. Thus, given $\mathsf{p}\in \mathcal{K}$
there exists $j^{\ast }\in \left\{ 1,...,N_{\epsilon }\right\} $ such that%
\begin{equation}
\left\Vert \mathsf{p}-\mathsf{c}_{j^{\ast }}\right\Vert _{\hat{\Phi}%
}<\epsilon \leq 1  \label{inaball}
\end{equation}%
as a consequence of (\ref{covering}). Writing then%
\begin{equation*}
\mathsf{p=c}_{j^{\ast }}+\left( \mathsf{p-c}_{j^{\ast }}\right)
\end{equation*}%
and using the fact that each $\Phi _{\mathsf{m}}$ is monotone increasing and
convex, we have successively%
\begin{eqnarray}
&&\Phi _{\mathsf{m}}\left( \left\vert p_{\mathsf{m}}\right\vert \right) 
\notag \\
&\leq &\frac{1}{2}\left( \Phi _{\mathsf{m}}\left( 2\left\vert \mathsf{c}%
_{j^{\ast },\mathsf{m}}\right\vert \right) +\Phi _{\mathsf{m}}\left(
2\left\vert p_{\mathsf{m}}-\mathsf{c}_{j^{\ast },\mathsf{m}}\right\vert
\right) \right)  \label{basicinequality} \\
&\leq &\frac{\eta }{2}\left( \Phi _{\mathsf{m}}\left( \left\vert \mathsf{c}%
_{j^{\ast },\mathsf{m}}\right\vert \right) +\Phi _{\mathsf{m}}\left(
\left\vert p_{\mathsf{m}}-\mathsf{c}_{j^{\ast },\mathsf{m}}\right\vert
\right) \right)  \notag
\end{eqnarray}%
for every $\left\vert \mathsf{m}\right\vert \geq \mathsf{M}$, where we used
the uniform global $\Delta _{2}$-condition in its original form (\ref%
{globalcondition}). Therefore, defining $\mathsf{M}_{\epsilon }^{\prime
}:=\max \left( \mathsf{M}_{\epsilon },\mathsf{M}\right) $ where $\mathsf{M}%
_{\epsilon }$ is given by (\ref{maximum}) and summing both sides of (\ref%
{basicinequality}) over $\left\vert \mathsf{m}\right\vert \geq \mathsf{M}%
_{\epsilon }^{\prime }$ we get%
\begin{eqnarray*}
&&\sum_{\left\vert \mathsf{m}\right\vert \geq \mathsf{M}_{\epsilon }^{\prime
}}\Phi _{\mathsf{m}}\left( \left\vert p_{\mathsf{m}}\right\vert \right) \\
&\leq &\frac{\eta }{2}\left( \sum_{\left\vert \mathsf{m}\right\vert \geq 
\mathsf{M}_{\epsilon }^{\prime }}\Phi _{\mathsf{m}}\left( \left\vert \mathsf{%
c}_{j^{\ast },\mathsf{m}}\right\vert \right) +\sum_{\left\vert \mathsf{m}%
\right\vert \geq \mathsf{M}_{\epsilon }^{\prime }}\Phi _{\mathsf{m}}\left(
\left\vert p_{\mathsf{m}}-\mathsf{c}_{j^{\ast },\mathsf{m}}\right\vert
\right) \right) \\
&\leq &\frac{\eta }{2}\left( \epsilon +\left\Vert \mathsf{p}-\mathsf{c}%
_{j^{\ast }}\right\Vert _{\hat{\Phi}}\right) \leq \eta \epsilon
\end{eqnarray*}%
as a consequence of (\ref{importantrelationbis}), (\ref{smalltail}) and (\ref%
{inaball}).

Let us now consider the case $\epsilon \in \left( 1,\infty \right) $, for
which we claim that an estimate similar to the preceding one holds true.
Indeed, we may repeat exactly the same argument as above with a finite
covering of $\mathcal{\bar{K}}$ of the form%
\begin{equation*}
\mathcal{\bar{K}\subseteq \cup }_{j=1}^{\hat{N}_{\epsilon }}\mathcal{B}%
\left( \epsilon ^{-1},\mathsf{\hat{c}}_{j}\right)
\end{equation*}%
with balls $\mathcal{B}\left( \epsilon ^{-1},\mathsf{\hat{c}}_{j}\right) $
of radius $\epsilon ^{-1}$ and centers $\mathsf{\hat{c}}_{j}\in \mathcal{%
\bar{K}}$, which leads to%
\begin{equation*}
\left\Vert \mathsf{p}-\mathsf{\hat{c}}_{\hat{\jmath}}\right\Vert _{\hat{\Phi}%
}<\epsilon ^{-1}\leq 1
\end{equation*}%
for some $\hat{\jmath}\in \left\{ 1,...,\hat{N}_{\epsilon }\right\} $.
Therefore, with the new decomposition%
\begin{equation*}
\mathsf{p=\hat{c}}_{\hat{\jmath}}+\left( \mathsf{p-\hat{c}}_{\hat{\jmath}%
\text{ }}\right)
\end{equation*}%
of every $\mathsf{p}\in \mathcal{K}$ we eventually obtain%
\begin{eqnarray*}
&&\sum_{\left\vert \mathsf{m}\right\vert \geq \mathsf{M}_{\epsilon }^{\prime
\prime }}\Phi _{\mathsf{m}}\left( \left\vert p_{\mathsf{m}}\right\vert
\right) \\
&\leq &\frac{\eta }{2}\left( \epsilon +\left\Vert \mathsf{p}-\mathsf{\hat{c}}%
_{\hat{\jmath}}\right\Vert _{\hat{\Phi}}\right) \leq \frac{\eta }{2}\left(
\epsilon +\epsilon ^{-1}\right) \leq \eta \epsilon
\end{eqnarray*}%
for some $\mathsf{M}_{\epsilon }^{\prime \prime }\in \mathbb{N}$. With $%
\mathsf{M}_{\epsilon }^{\prime \prime \prime }:=\max \left( \mathsf{M}%
_{\epsilon }^{\prime },\mathsf{M}_{\epsilon }^{\prime \prime }\right) $ we
therefore conclude that%
\begin{equation*}
\sum_{\left\vert \mathsf{m}\right\vert \geq \mathsf{M}_{\epsilon }^{\prime
\prime \prime }}\Phi _{\mathsf{m}}\left( \left\vert p_{\mathsf{m}%
}\right\vert \right) \leq \eta \epsilon
\end{equation*}%
for all $\mathsf{p}\in \mathcal{K}$ and each $\epsilon \in \left( 0,\infty
\right) $, which is (\ref{conditionquarto}). \ \ $\blacksquare $

\bigskip

\textsc{Remark.} Since the Orlicz functions $\Phi _{\mathsf{m}}\left(
t\right) =t^{\alpha _{\mathsf{m}}}$ with $\alpha _{\mathsf{m}}\in \left[
1,\infty \right) $ satisfy the uniform global $\Delta _{2}$-condition when
the sequence $\left( \alpha _{\mathsf{m}}\right) $ is bounded, Corollary 2
allows one to conclude that a bounded set $\mathcal{K}$ in Nakano's space $%
\mathsf{s}_{\hat{\Phi}}=l_{\hat{\Phi}}$ is relatively compact if, and only
if, the relation%
\begin{equation*}
\lim_{\mathsf{M}\rightarrow \infty }\sup_{\mathsf{p}\in \mathcal{K}%
}\sum_{\left\vert \mathsf{m}\right\vert \geq \mathsf{M}}\left\vert p_{%
\mathsf{m}}\right\vert ^{\alpha _{\mathsf{m}}}=0
\end{equation*}%
holds. In particular, if $\alpha _{\mathsf{m}}$ is independent of $\mathsf{m}
$, we retrieve the well-known fact that a bounded set $\mathcal{K}\subset l_{%
\mathbb{C}}^{\alpha }$ with $\alpha \in \left[ 1,\infty \right) $ is
relatively compact if, and only if, the condition%
\begin{equation*}
\lim_{\mathsf{M}\rightarrow \infty }\sup_{\mathsf{p}\in \mathcal{K}%
}\sum_{\left\vert \mathsf{m}\right\vert \geq \mathsf{M}}\left\vert p_{%
\mathsf{m}}\right\vert ^{\alpha }=0
\end{equation*}%
holds (see, e.g., Chapter 1 in \cite{Diestel}). We refer the reader to \cite%
{banomaak} for other criteria regarding the relative compactness of bounded
sets in modular sequence spaces.

\bigskip

We can modify the requirement (\ref{condition}) of Theorem 6 in order to
eliminate altogether the role of the $\Delta _{2}$-condition and thereby
allow for exponential growth of the functions in the sequence $\left( \Phi _{%
\mathsf{m}}\right) $. This is achieved by introducing an additional
parameter $\gamma \in \left( 0,\infty \right) $ that establishes a link
between the constants $\kappa $ and $\epsilon $ of Lemmas 1 and 2:

\bigskip

\textbf{Theorem 6 bis. }\textit{Let }$\left( \Phi _{\mathsf{m}}\right) $%
\textit{\ and }$\left( \Psi _{\mathsf{m}}\right) $\textit{\ be two sequences
of Orlicz functions, and let us assume that }%
\begin{equation}
\lim_{\left\vert \mathsf{m}\right\vert \rightarrow \infty
}\sup_{0<t\leqslant \Psi _{\mathsf{m}}^{-1}\left( 1\right) }\frac{\Phi _{%
\mathsf{m}}\left( \gamma t\right) }{\Psi _{\mathsf{m}}\left( t\right) }=0
\label{conditionterco}
\end{equation}%
\textit{for every }$\gamma \in \left( 0,\infty \right) $\textit{. Then there
exists the compact embedding}%
\begin{equation*}
\mathsf{s}_{\hat{\Psi}}\hookrightarrow \mathsf{s}_{\hat{\Phi}}.
\end{equation*}

\bigskip

\textbf{Proof. }According to the second part of the proof of Theorem 6, it
is sufficient to show that there exists $\mathsf{M}_{\epsilon ,\kappa }\in 
\mathbb{N}$ such that%
\begin{equation}
\displaystyle\sum\limits_{\left\vert \mathsf{m}\right\vert \geq \mathsf{M}_{\epsilon
,\kappa }}\Phi _{\mathsf{m}}\left( \frac{2\left\vert p_{\mathsf{m}%
}\right\vert }{\epsilon }\right) \leq 1  \label{relation}
\end{equation}%
for every $\epsilon \in \left( 0,\infty \right) $. By hypothesis, for every
finite $\mu >0$ and $\gamma >0$ there exists $\mathsf{M}_{\mu ,\gamma }\in 
\mathbb{N}$ such that%
\begin{equation*}
\Phi _{\mathsf{m}}\left( \frac{\gamma \left\vert p_{\mathsf{m}}\right\vert }{%
\kappa }\right) \leq \mu \Psi _{\mathsf{m}}\left( \frac{\left\vert p_{%
\mathsf{m}}\right\vert }{\kappa }\right)
\end{equation*}%
whenever $\left\vert \mathsf{m}\right\vert \geq \mathsf{M}_{\mu ,\gamma }$,
uniformly in $\mathsf{p}\in \mathcal{K}.$ Then, considering first $\epsilon
\in \left( 0,1\right] $ and choosing $\mu =\epsilon $, $\gamma =2\kappa
\epsilon ^{-1}$ in the preceding inequality we obtain%
\begin{equation*}
\sum_{\left\vert \mathsf{m}\right\vert \geq \mathsf{M}_{1,\epsilon ,\kappa
}}\Phi _{\mathsf{m}}\left( \frac{2\left\vert p_{\mathsf{m}}\right\vert }{%
\epsilon }\right) \leq \epsilon \sum_{\left\vert \mathsf{m}\right\vert \geq 
\mathsf{M}_{1,\epsilon ,\kappa }}\Psi _{\mathsf{m}}\left( \frac{\left\vert
p_{\mathsf{m}}\right\vert }{\kappa }\right) \leq \epsilon \leq 1
\end{equation*}%
for some $\mathsf{M}_{1,\epsilon ,\kappa }\in \mathbb{N}$, by invoking (\ref%
{relationter}). By the same token, for any $\epsilon \in \left( 1,\infty
\right) $ and choosing $\mu =\epsilon ^{-1}$ with the same value of $\gamma $
as above we get%
\begin{equation*}
\sum_{\left\vert \mathsf{m}\right\vert \geq \mathsf{M}_{2,\epsilon ,\kappa
}}\Phi _{\mathsf{m}}\left( \frac{2\left\vert p_{\mathsf{m}}\right\vert }{%
\epsilon }\right) \leq \epsilon ^{-1}\leq 1
\end{equation*}%
for some $\mathsf{M}_{2,\epsilon ,\kappa }\in \mathbb{N}$. Relation (\ref%
{relation}) then follows with $\mathsf{M}_{\epsilon ,\kappa }:=\allowbreak\max \left( 
\mathsf{M}_{1,\epsilon ,\kappa },\mathsf{M}_{2,\epsilon ,\kappa }\right) $.
\ \ $\blacksquare $

\bigskip

\textsc{Remarks.} In the preceding section we observed that there is no
Schauder basis in $l_{\hat{\Phi}}$ when $\mathsf{s}_{\hat{\Phi}}\subset l_{%
\hat{\Phi}}$, so that the preceding arguments may not be applied in this
case. It is therefore unclear whether there exist compact embeddings of the
form $l_{\hat{\Psi}}\hookrightarrow l_{\hat{\Phi}}$ in general.

\bigskip

Let us now illustrate Theorems 6 and 6 bis by means of two examples.

\bigskip

\textsc{Example 1.} We first consider the sequences $\left( \Phi _{\mathsf{m}%
}\right) $ and $\left( \Psi _{\mathsf{m}}\right) $ given by%
\begin{equation}
\Phi _{\mathsf{m}}\left( t\right) =\frac{t^{\alpha }}{\left\vert \mathsf{m}%
\right\vert +1}  \label{choice}
\end{equation}%
and%
\begin{equation}
\Psi _{\mathsf{m}}\left( t\right) =\exp \left[ t^{\beta }\right] -1
\label{choicebis}
\end{equation}%
respectively, with $\alpha ,\beta \in \left[ 1,\infty \right) $ and $\alpha
\geq \beta .$ Then we have%
\begin{equation*}
\sup_{0<t\leqslant \Psi _{\mathsf{m}}^{-1}\left( 1\right) }\frac{\Phi _{%
\mathsf{m}}\left( t\right) }{\Psi _{\mathsf{m}}\left( t\right) }=\frac{1}{%
\left\vert \mathsf{m}\right\vert +1}\sup_{0<t\leqslant \ln ^{\frac{1}{\beta }%
}(2)}\frac{t^{\alpha }}{\exp \left[ t^{\beta }\right] -1}<\infty
\end{equation*}%
since the ratio of these two functions may be extended by continuity at $%
t=0_{+}$. Therefore (\ref{condition}) holds, so that the embedding%
\begin{equation}
\mathsf{s}_{\hat{\Psi}}\hookrightarrow \mathsf{s}_{\hat{\Phi}}
\label{embedding5}
\end{equation}%
is compact. Under the same conditions regarding $\alpha $ and $\beta $ we
get an identical result if we replace (\ref{choice}) by 
\begin{equation*}
\Phi _{\mathsf{m}}\left( t\right) =\frac{\exp \left[ t^{\alpha }\right] -1}{%
\left\vert \mathsf{m}\right\vert +1},
\end{equation*}%
for (\ref{conditionterco}) holds.

\bigskip

\textsc{Example 2.} We obtain a similar result if%
\begin{equation*}
\Phi _{\mathsf{m}}\left( t\right) =\mu _{\left\vert \mathsf{m}\right\vert
}t^{\alpha }
\end{equation*}%
and%
\begin{equation*}
\Psi _{\mathsf{m}}\left( t\right) =\nu _{\left\vert \mathsf{m}\right\vert
}t^{\beta }
\end{equation*}%
when $\alpha ,\beta \in \left[ 1,\infty \right) $ and $\alpha \geq \beta $,
where $\left( \mu _{\left\vert \mathsf{m}\right\vert }\right) $ and $(\nu
_{\left\vert \mathsf{m}\right\vert })$ are two numerical sequences of
positive numbers satisfying%
\begin{equation*}
\lim_{\left\vert \mathsf{m}\right\vert \rightarrow \infty }\mu _{\left\vert 
\mathsf{m}\right\vert }\nu _{\left\vert \mathsf{m}\right\vert }^{-\frac{%
\alpha }{\beta }}=0.
\end{equation*}%
Indeed, in this case we have 
\begin{equation*}
\sup_{0<t\leqslant \Psi _{\mathsf{m}}^{-1}\left( 1\right) }\frac{\Phi _{%
\mathsf{m}}\left( t\right) }{\Psi _{\mathsf{m}}\left( t\right) }=\mu
_{\left\vert \mathsf{m}\right\vert }\nu _{\left\vert \mathsf{m}\right\vert
}^{-\frac{\alpha }{\beta }}.
\end{equation*}%
For instance, 
\begin{equation*}
\Phi _{\mathsf{m}}\left( t\right) =t^{2}
\end{equation*}%
and%
\begin{equation*}
\Psi _{\mathsf{m}}\left( t\right) =\left( 1+\mathsf{m}^{2}\right) t^{2}
\end{equation*}%
will do, from which we infer the known fact that the embedding%
\begin{equation*}
h_{\mathbb{C}}^{1}\hookrightarrow l_{\mathbb{C}}^{2}
\end{equation*}%
is compact, where $h_{\mathbb{C}}^{1}$ refers to the space introduced in the
remark following the proof of Theorem 6. This could also have been inferred
from (\ref{compactembeddingter}) by virtue of the composition%
\begin{equation*}
h_{\mathbb{C}}^{1}\hookrightarrow l_{\mathbb{C}}^{1}\rightarrow l_{\mathbb{C}%
}^{2}
\end{equation*}%
where the second embedding is continuous.

\bigskip

Theorem 6 leads to further results concerning the modular classes of Theorem
5, with the appropriate restrictions:

\bigskip

\textbf{Corollary 3. }\textit{Let us consider the sequence }$\left( \Phi _{%
\mathsf{m}}\right) $\textit{\ defined by }%
\begin{equation}
\Phi _{\mathsf{m}}(t):=\mu _{\left\vert \mathsf{m}\right\vert }\Phi \left(
t\right)  \label{definition}
\end{equation}%
\textit{where }$\Phi $\textit{\ is a fixed Orlicz function satisfying the
global }$\Delta _{2}$\textit{-condition}, \textit{and} $\left( \mu
_{\left\vert \mathsf{m}\right\vert }\right) $\textit{\ a numerical sequence
satisfying }$\mu _{\left\vert \mathsf{m}\right\vert }>0$\textit{\ for every }%
$\mathsf{m}$\textit{.} \textit{Let }$\left( \Psi _{\mathsf{m}}\right) $ 
\textit{be a sequence of} \textit{Orlicz functions and let us assume that
the inequality}%
\begin{equation}
\Phi \left( t\right) \leq c\Psi _{\mathsf{m}}\left( t\right)
\label{inequalitysetimo}
\end{equation}%
\textit{holds for some finite }$c>0$,\textit{\ each }$t\in \left( 0,\tau %
\right] $ \textit{where }$\tau $ \textit{is given by (\ref{boundedness}),
and every sufficiently large }$\left\vert \mathsf{m}\right\vert $\textit{.
Then the following statements hold:}

\textit{(a) If}%
\begin{equation}
\lim_{\left\vert \mathsf{m}\right\vert \rightarrow \infty }\mu _{\left\vert 
\mathsf{m}\right\vert }=0,  \label{limit}
\end{equation}%
\textit{then} \textit{there exists the compact embedding}%
\begin{equation}
\mathsf{s}_{\hat{\Psi}}\hookrightarrow \mathsf{s}_{\hat{\Phi}}.
\label{compactembeddingquarto}
\end{equation}

\textit{(b) Conversely, let us assume in addition that each right-derivative 
}$\Psi _{r,\mathsf{m}}$ \textit{satisfies }$\Psi _{r,\mathsf{m}}\left(
0\right) =0$ \textit{and }$\lim_{t\rightarrow \infty }\Psi _{r,\mathsf{m}r,%
\mathsf{m}}\left( t\right) =\infty $, \textit{and that}%
\begin{equation}
\liminf_{\left\vert \mathsf{m}\right\vert \rightarrow \infty }\Psi _{\mathsf{%
m}}^{\ast }\left( t\right) >0  \label{constraintter}
\end{equation}%
\textit{for every} $t>0$. \textit{If there exists the compact embedding (\ref%
{compactembeddingquarto}), then (\ref{limit}) holds.}

\bigskip

\textbf{Proof.} The uniform global $\Delta _{2}$\textit{-}condition imposed
on $\Phi $ implies that the same condition holds for each $\Phi _{\mathsf{m}%
} $. Moreover we have%
\begin{equation*}
\frac{\Phi _{\mathsf{m}}\left( t\right) }{\Psi _{\mathsf{m}}\left( t\right) }%
\leq c\mu _{\left\vert \mathsf{m}\right\vert }
\end{equation*}%
uniformly in $t\in \left( 0,\tau \right] $ as a consequence of (\ref%
{inequalitysetimo}), so that%
\begin{equation*}
\lim_{\left\vert \mathsf{m}\right\vert \rightarrow \infty
}\sup_{0<t\leqslant \Psi _{\mathsf{m}}^{-1}\left( 1\right) }\frac{\Phi _{%
\mathsf{m}}\left( t\right) }{\Psi _{\mathsf{m}}\left( t\right) }\leq
\lim_{\left\vert \mathsf{m}\right\vert \rightarrow \infty
}\sup_{0<t\leqslant \tau }\frac{\Phi _{\mathsf{m}}\left( t\right) }{\Psi _{%
\mathsf{m}}\left( t\right) }=0
\end{equation*}%
if (\ref{limit}) holds. Therefore, there exists the compact embedding $%
\mathsf{s}_{\hat{\Psi}}\hookrightarrow \mathsf{s}_{\hat{\Phi}}$ according to
Theorem 6, which proves Statement (a).

In order to prove Statement(b) we assume that (\ref{limit}) does not hold
and show that a compact embedding of the form (\ref{compactembeddingquarto})
is impossible. Under this hypothesis the sequence $\left( \mu _{\left\vert 
\mathsf{m}\right\vert }\right) $ may be either unbounded or bounded. In the
first case, there may or may not exist a continuous embedding of $\mathsf{s}%
_{\hat{\Psi}}$ into $\mathsf{s}_{\hat{\Phi}}$ according to Proposition 2 and
Remark (2) following it. If such an embedding does not exist, then \textit{a
fortiori }a compact embedding does not exist either. Now the conditions of
Statement (b) ensure that the sequence $\left( \mathsf{e}_{\mathsf{m}%
}\right) $ of basis vectors converges weakly to zero in $\mathsf{s}_{\hat{%
\Psi}}$ according to Proposition 1 of the preceding section. Therefore, if
there exists a continuous embedding of $\mathsf{s}_{\hat{\Psi}}$ into $%
\mathsf{s}_{\hat{\Phi}}$, which is also the case if $\left( \mu _{\left\vert 
\mathsf{m}\right\vert }\right) $ remains bounded according to Theorem 5, its
non-compactness will follow if we can show that $\left( \mathsf{e}_{\mathsf{m%
}}\right) $\textsf{\ does} not converge strongly to zero in the target space 
$\mathsf{s}_{\hat{\Phi}}$. But from (\ref{norm}), the fact that $\left( 
\mathsf{e}_{\mathsf{m}}\right) _{\mathsf{n}}=\delta _{\mathsf{m,n}}$ for all 
$\mathsf{m,n}\in \mathbb{Z}$ and (\ref{definition}) we obtain%
\begin{equation*}
\left\Vert \mathsf{e}_{\mathsf{m}}\right\Vert _{\hat{\Phi}}=\inf \left\{
\rho >0:\rho \geq \frac{1}{\Phi _{\mathsf{m}}^{-1}\left( 1\right) }\right\} =%
\frac{1}{\Phi _{\mathsf{m}}^{-1}\left( 1\right) }=\frac{1}{\Phi ^{-1}\left(
\mu _{\left\vert \mathsf{m}\right\vert }^{-1}\right) }.
\end{equation*}%
Therefore, the behavior of $\Phi ^{-1}$ for large arguments as the inverse
of an Orlicz function implies that the preceding expression does not go to
zero as $\left\vert \mathsf{m}\right\vert \rightarrow \infty $ when (\ref%
{limit}) does not hold. \ \ $\blacksquare $

\bigskip

If the sequence $\left( \mu _{\left\vert \mathsf{m}\right\vert }\right) $
depends on an additional scalar parameter, then under some conditions
Theorem 6 implies the existence of a compact embedding involving spaces of
Orlicz-Sobolev type. More specifically, the following result holds:

\bigskip

\textbf{Corollary 4.}\textit{\ Let }$\left( \Phi _{k,\mathsf{m}}\right) $%
\textit{\ be a sequence of Orlicz functions given by}%
\begin{equation*}
\Phi _{k,\mathsf{m}}\left( t\right) =\mu _{k,\left\vert \mathsf{m}%
\right\vert }\Phi \left( t\right)
\end{equation*}%
\textit{where }$\Phi $\textit{\ is a fixed Orlicz function satisfying the
global }$\Delta _{2}$\textit{-condition}.\ \textit{In the preceding
expression}, $\left( \mu _{k,\left\vert \mathsf{m}\right\vert }\right) $%
\textit{\ is a numerical sequence defined for each }$k\in \mathbb{R}$ 
\textit{satisfying }$\mu _{k,\left\vert \mathsf{m}\right\vert }>0$ \textit{%
for every }$\mathsf{m}$. \textit{Furthermore, let us assume that the
functional relation}%
\begin{equation}
\mu _{k,\left\vert \mathsf{m}\right\vert }=\mu _{k-k^{\prime },\left\vert 
\mathsf{m}\right\vert }\times \mu _{k^{\prime },\left\vert \mathsf{m}%
\right\vert }  \label{functionalrelation}
\end{equation}%
\textit{holds for all }$k,k^{\prime }\in \mathbb{R}$\textit{, and that }%
\begin{equation}
\lim_{\left\vert \mathsf{m}\right\vert \rightarrow \infty }\mu _{k-k^{\prime
},\left\vert \mathsf{m}\right\vert }=0  \label{limitter}
\end{equation}%
\textit{whenever }$k^{\prime }>k$, \textit{uniformly in these variables.\
Then the embedding} 
\begin{equation}
\mathsf{s}_{k^{\prime },\hat{\Phi}}\hookrightarrow \mathsf{s}_{k,\hat{\Phi}}
\label{compactembedding}
\end{equation}%
\textit{is compact, where }$\mathsf{s}_{k,\hat{\Phi}}$\textit{\ is the
modular space generated by the sequence} $\left( \Phi _{k,\mathsf{m}}\right) 
$.

\bigskip

\textbf{Proof.} We proceed as in the proof of Corollary 3 to get%
\begin{equation*}
\frac{\Phi _{k,\mathsf{m}}\left( t\right) }{\Phi _{k^{\prime },\mathsf{m}%
}\left( t\right) }=\mu _{k-k^{\prime },\left\vert \mathsf{m}\right\vert }
\end{equation*}%
for each $\mathsf{m}$ uniformly in $t$ by virtue of (\ref{functionalrelation}%
), so that (\ref{limitter}) implies%
\begin{equation*}
\lim_{\left\vert \mathsf{m}\right\vert \rightarrow \infty
}\sup_{0<t\leqslant \Phi _{k^{\prime },\mathsf{m}}^{-1}\left( 1\right) }%
\frac{\Phi _{k,\mathsf{m}}\left( t\right) }{\Phi _{k^{\prime },\mathsf{m}%
}\left( t\right) }=0.
\end{equation*}%
Theorem 6 then entails the existence of (\ref{compactembedding}). \ \ $%
\blacksquare $

\bigskip

\textsc{Remark.} The relation between the corollary and the statements about
compact embeddings for spaces of Orlicz-Sobolev type we alluded to above is
apparent if we choose%
\begin{equation*}
\mu _{k,\left\vert \mathsf{m}\right\vert }=\left( 1+\Phi \left( \left\vert 
\mathsf{m}\right\vert \right) \right) ^{k},
\end{equation*}%
which satisfies all the required hypotheses. Thus, in this case we have 
\begin{equation*}
\Phi _{k,\mathsf{m}}\left( t\right) =\left( 1+\Phi \left( \left\vert \mathsf{%
m}\right\vert \right) \right) ^{k}\Phi \left( t\right) .
\end{equation*}%
It is to be noted that compactness in (\ref{compactembedding}) only emerges
when the so-called generalized order of differentiability $k^{\prime }$
associated with $\mathsf{s}_{k^{\prime },\hat{\Phi}}$ is strictly larger
than that of the target space $\mathsf{s}_{k,\hat{\Phi}}$. We refer the
reader to \cite{Vuillermotbis} and to the remark at the very end of Section
2 for more details regarding this terminology. There are of course plenty of
other possible choices for $\left( \mu _{k,\left\vert \mathsf{m}\right\vert
}\right) $, for instance $\mu _{k,\left\vert \mathsf{m}\right\vert }=\exp %
\left[ k\left\vert \mathsf{m}\right\vert \right] $.

\bigskip

\textbf{Acknowledgments.} The author is indebted to M. Mastylo for useful
comments regarding the topics treated in this article. He also thanks the
Funda\c{c}\~{a}o para a Ci\^{e}ncia e a Tecnologia of the Portugu\^{e}s
Government (FCT) for its financial support under grant UID/04561/2025.

\bigskip

Data availability: My manuscript has no associated data.

\bigskip

\textbf{Declarations}:

Ethical approval: Not applicable.

Conflict of interest: The author declares no competing interest.


\begin{thebibliography}{99}
\bibitem{banomaak} \textsc{Bandaliyev, R. A., Omarova, K. K., Akhundova, S.
M.}, \textit{Compactness criterion in weighted Musielak-Orlicz sequence
spaces and applications}, Transactions of A. Razmadze Mathematical Institute 
\textbf{180}, 1, 13-24 (2026).

\bibitem{Bisgard} \textsc{Bisgard, J., }\textit{A compact embedding for
sequence spaces,} Missouri Journal of Mathematical Sciences \textbf{24}, 2,
182-189 (2012).

\bibitem{Diharu} \textsc{Diening, L., Harjulehto, P., H\"{a}st\"{o}, P., R%
\r{u}\v{z}i\v{c}ka, M.,} \textit{Lebesgue and Sobolev spaces with variable
exponents}, Lecture Notes in Mathematics \textbf{2017}, Springer, New York
(2011).

\bibitem{Diestel} \textsc{Diestel, J.,} \textit{Sequences and Series in
Banach Spaces,} Graduate Texts in Mathematics \textbf{92}, Springer, New
York (1984).

\bibitem{Edmundsevans} \textsc{Edmunds, D. E., Evans, W. D.,} \textit{%
Spectral Theory and Differential Operators,} Oxford Mathematical Monographs,
Clarendon Press (1987).

\bibitem{Krasnoselskiirut} \textsc{Krasnosel'skii, M. A., Rutickii, Ya. B.,} 
\textit{Convex Functions and Orlicz Spaces}, P. Noordhoff, Ltd. Groningen
(1961).

\bibitem{Langnek} \textsc{Lang, J., Nekvinda, A.,} \textit{Embeddings
between sequence variable Lebesgue spaces, strict and finitely strict
singularity}, Mathematische Nachrichten, \textbf{298}, 9, 2926-2941 (2025).

\bibitem{Lindenstrausstza} \textsc{Lindenstrauss, J., Tzafriri, L.,} \textit{%
Classical Banach Spaces I: Sequence Spaces,} Ergebnisse der Mathematik und
ihrer Grenzgebiete \textbf{92}, Springer, New York (1977).

\bibitem{Musielak} \textsc{Musielak, J.,} \textit{Orlicz Spaces and Modular
Spaces,} Lecture Notes in Mathematics \textbf{1034}, Springer, New York
(1983).

\bibitem{Nakano} \textsc{Nakano, H.,}\textit{\ Modulared sequence spaces},
Proceedings of the Japan Academy, \textbf{27}, 9, 508-512 (1951).

\bibitem{Raoren} \textsc{Rao, M. M., Ren, Z. D.,} \textit{Theory of Orlicz
Spaces,} Pure and Applied Mathematical Series \textbf{146}, Marcel Dekker,
Inc., New York (1991).

\bibitem{Raorenbis} \textsc{Rao, M. M., Ren, Z. D.,}\textit{\ Applications
of Orlicz Spaces,} Pure and Applied Mathematical Series \textbf{250}, Marcel
Dekker, Inc., New York (2002).

\bibitem{Vuillermot} \textsc{Vuillermot, P.-A.,} \textit{On the relaxation
to equilibrium of a quantum oscillator interacting with a radiation field, }%
Communications in Mathematical Sciences \textbf{23}, 6, 1519-1546 (2025).

\bibitem{Vuillermotbis} \textsc{Vuillermot, P.-A., }\textit{Compact
embeddings and Pitt's property for a class for weighted sequence spaces of
Sobolev type,} Communications in Mathematical Analysis and Applications 
\textbf{4}, 2, 285-295 (2025).

\bibitem{Vuillermotter} \textsc{Vuillermot, P.-A., }\textit{On the
compactness of embeddings for a class of weighted Orlicz-Sobolev sequence
spaces,} Journal of Mathematical Sciences,
\url{https://doi.org/10.1007/s10958-025-08123-2} (2025).
\end{thebibliography}
\end{document}